\documentclass[twoside, UTF8]{article}
\usepackage{tikz}
\usetikzlibrary{matrix}
\usepackage[all]{xy}
\usepackage{graphics}
\usepackage{subfigure}
\usepackage{amsmath,amsopn,amssymb,epsfig,psfrag}
\usepackage{color,xcolor}
\usepackage{pdflscape}
\usepackage{float} %±íͼ¸úÎÄ×Öºó

\title{Large-Scale Algebraic Riccati Equations with High-Rank Nonlinear Terms and Constant Terms}

\date{}
\author{
Bo Yu\thanks{School of Science, Hunan University of Technology, Zhuzhou 412007, China; {\tt  boyu\_hut@126.com; dongning\_158@sina.com}}
\and Ning Dong\footnotemark[1] }
\begin{document}
%
%-------------------------------------------------------------------------------
%
%abb1.tex, revised from abb.tex; December 10, 1990.
%
%\newcommand{\bs}{\begin{slide}{}}
%\newcommand{\es}{\end{slide}} \newcommand{\be}{\begin{eqnarray*}}
%\newcommand{\ee}{\end{eqnarray*}}
%\newcommand{\be}[1]{\begin{eqnarray} \label{#1}}
%\newcommand{\ee}{\end{eqnarray}}
%\newcommand{\bt}{\begin{tabular}}
%\newcommand{\et}{\end{tabular}} \newcounter{bean}
%\newcommand{\bl}[1]{\begin{list}{#1}{\usecounter{bean}}} \newcommand{\el}{\end{list}}

\newcommand{\bel}[1]{\begin{equation}\label{#1}} \newcommand{\eel}{\end{equation}}

\newcommand{\sectionbreak}{\clearpage}

\def\xa{x_{1}}
\def\xb{x_{2}} \def\xc{x_{3}} %\newcommand{\ebs}{\end{slide}\begin{slide}{}}
\newtheorem{lem}{Lemma}[section]
\newtheorem{theorem}[lem]{Theorem}
\newtheorem{definition}[lem]{Definition}
\newtheorem{corollary}[lem]{Corollary}
\newtheorem{example}{Example}[section]
\newtheorem{remark}{Remark}[section]
\def\sr{{\cal R}} \def\cs{{\cal C}} \def\vx{{\bf x}}
%-------------------------------------------------------------------------------
%
\maketitle

\begin{abstract}
For large-scale discrete-time algebraic Riccati equations (DAREs) with high-rank nonlinear and constant terms, the stabilizing solutions are no longer numerically low-rank, resulting in the obstacle in the computation and storage. However, in some proper control problems such as power systems, the potential structure of the state matrix -- banded-plus-low-rank, might make the large-scale computation essentially workable. In this paper, a factorized structure-preserving doubling algorithm (FSDA) is developed under the frame of the banded inverse of nonlinear and constant terms. The detailed iterations format, as well as a deflation process of FSDA, are analyzed in detail. A technique of partial truncation and compression is introduced to shrink the dimension of columns of low-rank factors as much as possible. The computation of residual, together with the termination condition of the structured version, is also redesigned.
\end{abstract}



%\pagestyle{myheadings} \thispagestyle{plain}
%\markboth{B. Yu,  \& E.K.-W. Chu}{Large-Scale Algebraic Riccati Equations with High-Rank Constant Terms}

\section{Introduction}

Consider the LTI control system in discrete-time
\[
x(t+1) = Ax(t) + Bu(t), \ \ \ y(t) = Cx(t)
\]
with $A \in \mathbb{R}^{N \times N}$, $B \in \mathbb{R}^{N \times m}$ and $C \in \mathbb{R}^{l \times N}$ and $m,l \leq N$. The linear quadratic regulator (LQR) control minimizes
\[
J(x,u) = \sum_{t=0}^\infty \left[ x(t)^\top H x(t) + u(t)^\top R u(t) \right].
\]
The corresponding optimal control $u(t) = F x(t)$ and the feedback gain $F=-(R + B^\top X B)^{-1} B^\top X A$ can be expressed in terms of the unique positive semidefinite stabilzing solution $X$ of the descrete-time algebraic Riccati equation (DARE) \cite{af65, cflw,lr,m}:
\begin{eqnarray}\label{d}
\mathcal{D}(X) &=& -X + A^\top X (I + GX)^{-1} A + H.  \end{eqnarray}

The stabilizing solution of DARE \eqref{d} is of great importance in optimal control and has been an active area of research for the last several decades. Dozens of methods including the classical ones as well as the state-of-the-art ones have been devised to solve the equation in a numerically stable manner. See \cite{bf98, cfl, cflw, dld08, dld10, dz14, lr, l79, lc14, m, k68, zl20} and references therein for more details.

In many large-scale control problems, the matrix $G = BR^{-1} B^\top$ in the nonlinear term and $H=C^\top T^{-1} C$ in the constant term are of low-rank structure with $B \in \mathbb{R}^{N \times m^h}$, $R \in \mathbb{R}^{m^h \times m^h}$, $C \in \mathbb{R}^{m^g \times N}$, $T \in \mathbb{R}^{m^g \times m^g}$ and $m^g, m^h \ll N$. Then the unique positive definite stabilizing
solution in the DARE \eqref{d} or its dual equation can be approximated numerically by a low-rank matrix \cite{bs, cw15}. When DARE \eqref{d} has  a {\it high-rank} constant term $H$, the stabilizing solution is no longer numerically low-rank and its storage and outputting are nontrivial. By noting the remained low-rank structure of $G$, an adapted version of the doubling algorithm, i.e. SDA\_h was presented to solve the large-scale DARE efficiently \cite{yfc19}. The main idea behind that is to make full use of the numerical low-rank of the stabilizing solution in the dual equation to estimate the residual of the original DARE, so that the SDA\_h can realize the evaluation of the residual and the output of the feedback gain in a low-rank way. An interesting question up to now might be that
\begin{itemize}
  \item[] can SDA solve the large-scale DAREs efficiently when both $G$ and $H$ are of high-rank?
\end{itemize}

The main difficulty, in this case, lies in that the stabilizing solutions both in DARE \eqref{d} and its dual are not of low-rank structure, making the direct application of SDA\_h difficult for large-scale problems, let alone the estimation of DARE's residual and the realization of algorithmic termination. This paper attempts to cope with this obstacle to some extent. Rather than answer the above question entirely, DARE \eqref{d} with the banded-plus-low-rank structure
\begin{equation}\label{sa}
A = D^A + L^A_{10} K^A (L^A_{20})^\top
\end{equation}
is concerned, where $D^A\in \mathbb{R}^{N \times N}$ is the banded matrix, $L^A_{10}$, $L^A_{20}\in \mathbb{R}^{N \times m^a}$ are low-rank matrices and $K^A\in \mathbb{R}^{m^a \times m^a}$ is the kernel matrix with $ m^a\ll N$. Such a structure has some proper applications in the power system \cite{FMV11, MLP96, RM06}. Besides, the high-rank nonlinear item and the constant item are assumed to be
\begin{equation}\label{sgh}
G = D^G +L^GK^G(L^G)^\top, \ \ H = D^H +L^HK^H(L^H)^\top,
\end{equation}
where $D^G$, $D^H\in \mathbb{R}^{N \times N}$ are nonsingular
banded matrices, $L^G\in \mathbb{R}^{N \times m^g}$, $L^H\in \mathbb{R}^{N \times m^h}$, $K^G\in \mathbb{R}^{m^g \times m^g}$, $K^H\in \mathbb{R}^{m^h \times m^h}$ and $m^g, m^h\ll N$.
Moreover, from the viewpoint of applications, banded matrices $D^A$, $D^G$, $D^H$ are generally of a banded inverse (see %including kinds of matrices such as wavelet matrix \cite{SN96}, digital filter matrix \cite{k85}, CMV matrix\cite{cmv03} and Green's matrix\cite{OZS10} etc. See
\cite{dm93, cmv03, k85, km00, S10, S11, SN96, OZS10} as well as their references).

The main contributions in this paper are from the following aspects:
\begin{itemize}
  \item It is first to develop SDA to the factorized form --- FSDA --- to solve the large-scale DAREs with both high-rank $H$ and $G$, where the stabilizing solutions in DARE and its dual are no longer numerically low-rank.

\item The structure of the FSDA iterative sequence is explicitly revealed to consist of two parts --- the banded part and the low-rank part. The banded part can iterate independently while the low-rank part relies heavily on the product of the banded part and the low-rank part.

\item A deflation process of the low-rank factors is proposed to reduce the column number of the low-rank part. The conventional truncation and compression in \cite{cw15, yfc19} for the whole low-rank factor seems not to work as it destroys the implicit structure and makes the subsequent deflation infeasible. Instead, a partial truncation and compression (PTC) technique is then devised to impose merely on the exponentially increasing part (after deflation), effectively slimming the dimension of the columns in low-rank factors.

\item The termination of FSDA is designed to comprise two parts. The relatively easy-computing residual of the banded part makes up the pre-termination, followed by the actual termination condition from the residual of the low-rank factors. In this way, the latter time-consuming condition could be reduced as much as possible, cutting down the whole computational complexity of FSDA.
\end{itemize}

The whole paper is organized as follows. Section 2 describes the iteration format of FSDA for DAREs \eqref{d} with high-rank nonlinear and constant terms. The deflation process for the low-rank factors and kernels is given in Section 3. Section 4 dwells on the technique of PTC to slim the dimension of the columns of low-rank factors and kernels. The way to compute the residual, as well as the concrete implementation of FSDA, is described in Section 5. Numerical experiments are listed in Section 6 to show the effectiveness of FSDA.

{\bf Notation.} $I_{{}_N}$ (or simply $I$) is the $N\times N$ identity matrix. For a matrix $A\in \mathbb{R}^{N\times N}$, $\rho(A)$ denotes the spectral radius of $A$.  For symmetric matrices $A$ and $B\in\mathbb{R}^{N\times N}$, we say $A>B$ ($A\geq B$) if $A-B$ is a positive definitive (positive semi-definite) matrix. Additionally, the Sherman-Morrison-Woodbury formula (SMWF) (see \cite{gv96} for example),
$ %\bel{smwf}
(M + UDV^\top)^{-1} = M^{-1} - M^{-1} U (D^{-1} + V^\top M^{-1} U)^{-1} V^\top M^{-1}
$
%\] %\eel
is required in the analysis of iterative format.

%%%%%%%%%%%%%%%%%%%%%%%%%%%%%%%%%%%%%%%%%%%%%%%%
\section{SDA and the structured iteration for DARE}
%%%%%%%%%%%%%%%%%%%%%%%%%%%%%%%%%%%%%%%%%%%%%%%
For DARE
\[\mathcal{D}(X)= -X + A^\top X (I + GX)^{-1} A + H =0 \] and its dual equation
\begin{equation}\label{dd}
\mathcal{D}_a (Y) = -Y + AY (I + HY)^{-1} A^\top + G = 0,
\end{equation}
%with   \textit{high-rank} nonlinear term $G$ and constant term $H$ in \eqref{sgh}.
%$$ G = D^G +L^GK^G(L^G)^\top \geq 0, \ \ H = D^H +L^HK^H(L^H) \geq 0, $$
%where $D^G$, $D^H\in \mathbb{R}^{N \times N}$ are banded matrices with banded inverse and $L^G\in \mathbb{R}^{N \times m^g}$, $L^H\in \mathbb{R}^{N \times m^h}$, $K^G\in \mathbb{R}^{m^g \times m^g}$, $K^H\in \mathbb{R}^{m^h \times m^h}$ with $m^g,m^h \ll N$.
%For large-scale DAREs, we further assume that $m^g,m^h \ll N$ and $A, M, N$ are BMBI with respective bandwidth $b^a$, $b^g$ and $b^h$.
%For the large-scale CARE \eqref{c} with $m \ll n$, we need the matrices being efficiently invertible and multiplicable to vectors; see Section~4. CAREs are mathematically equivalent to the DARE \eqref{d}, through the Cayley transformation \eqref{cayleya}.
SDA \cite{cflw} proceeds for $k\geq1$
\bel{sda}
\left\{
\begin{array}{rcl}
G_{k} &=& G_{k-1} + A_{k-1} (I+G_{k-1} H_{k-1} )^{-1}G_{k-1} A_{k-1}^\top, \\
H_{k} &=& H_{k-1} + A_{k-1}^\top H_{k-1}(I+G_{k-1} H_{k-1})^{-1} A_{k-1}, \\
A_{k} &=& A_{k-1} (I+G_{k-1} H_{k-1})^{-1} A_{k-1},
\end{array}
\right.
\eel
with $A_0=A$, $G_0=G$, $H_0=H$. %Under the initial structure of \eqref{sa} and \eqref{sgh},  iterative scheme of FSDA will be constructed as below.

\subsection{FSDA for high-rank terms}
Given banded matrices $D_0^A$, $D_0^G$ and $D_0^H$ and low-rank matrices $L_{10}^A$ and $L_{20}^A$ in structured initial matrices \eqref{sa} and \eqref{sgh}, the frame of the FSDA is described inductively as follows:
\bel{agh}
A_k = D_k^A + L_{1,k}^AK_k^A (L_{2,k}^A)^\top, \ \ \
G_k = D_k^G + L_k^G K_k^G (L_k^G)^\top, \ \ \ H_k = D_k^H + L_k^H K_k^H (L_k^H)^\top
\eel
with sparse banded matrices $D^A_k, D^G_k, D^H_k \in \mathbb{R}^{N\times N}$, low-rank factors $L^A_{1,k}\in \mathbb{R}^{N\times m^{a_1}_k}$, $L^A_{2,k}\in \mathbb{R}^{N\times m^{a_2}_k}$, $L^G_k\in \mathbb{R}^{N\times m^{g}_k}$, $L^H_k\in \mathbb{R}^{N\times m^{h}_k}$, kernel matrices $K^A_k\in \mathbb{R}^{m^{a_1}_k\times m^{a_2}_k}$, $K^G_k\in \mathbb{R}^{m^{g}_k\times m^{g}_k}$, $K^H_k\in \mathbb{R}^{m^{h}_k\times m^{h}_k}$ and $m^{a_1}_k, m^{a_2}_k, m^g_k, m^h_k\ll N$. Without loss of
generality, we assume that $m^{a_1}_0=m^{a_2}_0 \equiv m^{a}$ %the initial column number of $L^A_{10}$ is assumed to be the same as that of $L^A_{20}$
and $K^A_0=I_{m^a}$. Otherwise, a redefinition of $L^A_{20}:= K^A_0L^A_{20}$ and $K^A_0:= I_{m^a}$ can fulfill the assumption. Besides,  $K^G$ and $K^H$ are assumed to be zero for simplicity, and it will not destroy the structure.

We will elaborate the concrete format of banded matrices and low-rank factors for $k=1$ and $k\geq 2$. Note that banded parts are capable of iterating independently, regardless of low-rank parts and kernels.

\subsubsection*{Case for $k=1$. }
Insert the initial matrices $D_0^A$, $D_0^G$ and $D_0^H$ and low-rank matrices $L_{10}^A$ and $L_{20}^A$ into SDA \eqref{sda}. It follows from the SMW formula that banded parts are
\bel{sda1}
\left.
\begin{array}{rcl}
D^G_{1} &=& D^G_0 + D_0^{AGHG} (D^A_0)^\top, \\
D^H_{1} &=& D^H_0 + D_0^{A^{\top}HGH} D^A_0, \\
D^A_{1} &=& D_0^{AGH} D^A_0 = D^A_0 (D_0^{A^{\top}HG})^\top
\end{array}
\right.
\eel
with
\[ %\begin{equation}\label{daghahg}
\left.
\begin{array}{rcl}
D_0^{AGHG} =D_0^A(I_{{}_N}+D_0^GD_0^H)^{-1}D_0^G, && D_0^{A^{\top}HGH}= (D_0^A)^{\top}(I_{{}_N}+D_0^HD_0^G)^{-1}D_0^H,\\
 D_0^{AGH} =D_0^A(I_{{}_N}+D_0^GD_0^H)^{-1}, &&
 D_0^{A^{\top}HG}= (D_0^A)^{\top}(I_{{}_N}+D_0^HD_0^G)^{-1}.
\end{array}
\right.
%\end{equation}
\]
Low-rank factors in \eqref{agh} are
\begin{equation}\label{lagh1}
\left.
\begin{array}{rcl}
L_1^G =[L_{10}^A, \ D_0^{AGHG}L_{20}^A], & & L_1^H =[L_{20}^A, \ D_0^{A^{\top}HGH}L_{10}^A], \\
L_{1,1}^A =[L_{10}^A, D_0^{AGH}L_{10}^A], & & L_{2,1}^A =[L_{20}^A, D_0^{A^\top HG}L_{20}^A]
\end{array}
\right.
\end{equation}
and kernels in the low-rank parts are
\begin{eqnarray}
K_1^G &=& \left[
\begin{matrix}
(L_{20}^A)^\top D_0^{GHG}L_{20}^A & I_{m^{g}_0} \\
I_{m^{g}_0} & 0
\end{matrix} \right], \label{kg1}\\
K_1^H &=&\left[
\begin{matrix}
(L_{10}^A)^\top D_0^{HGH}L_{10}^A & I_{m^{h}_0} \\
I_{m^{h}_0} & 0
\end{matrix}\right], \label{kh1}\\
K_1^A &=&\left[
\begin{matrix}
(L_{20}^A)^\top D_0^{GH}L_{10}^A & I_{m^{g}_0} \\
I_{m^{h}_0} & 0
\end{matrix}\right]\label{ka1}
\end{eqnarray}
with
\[%\begin{equation}\label{dghhg}
\left.
\begin{array}{rcl}
D_0^{GHG} =(I_{{}_N}+D_0^GD_0^H)^{-1}D_0^G, \ \
D_0^{HGH} = (I_{{}_N}+D_0^HD_0^G)^{-1}D_0^H,\ \
 D_0^{GH} = (I_{{}_N}+D_0^GD_0^H)^{-1}
\end{array}
\right.
\]%\end{equation}
and $m^{g}_0=m^a$, $m^{h}_0=m^a$.

\subsubsection*{Case for general $k\geq 2$.}
Inserting banded matrices $D_{k-1}^{H}$, $D_{k-1}^{G}$ and $D_{k-1}^{A}$  and low-rank factors $L^G_{k-1}$, $L^H_{k-1}$, $L^A_{1,k-1}$ and $L^A_{2,k-1}$ into SDA \eqref{sda} again, banded matrices at the $k$-th iteration are
\bel{sdak}
\left.
\begin{array}{rcl}
D^G_{k} &=& D^G_{k-1} + D_{k-1}^{AGHG} (D^A_{k-1})^\top, \\
D^H_{k} &=& D^H_{k-1} + D_{k-1}^{A^{\top}HGH} D^A_{k-1}, \\
D^A_{k} &=& D_{k-1}^{AGH} D^A_{k-1} = D^A_{k-1} (D_{k-1}^{A^{\top}HG})^\top
\end{array}
\right.
\eel
with
\[ %\begin{equation}\label{daghahgk}
\left.
\begin{array}{rcl}
D_{k-1}^{AGHG} =D_{k-1}^A(I_{{}_N}+D_{k-1}^GD_{k-1}^H)^{-1}D_{k-1}^G, && D_{k-1}^{A^{\top}HGH}= (D_{k-1}^A)^{\top}(I_{{}_N}+D_{k-1}^HD_{k-1}^G)^{-1}D_{k-1}^H,\\
 D_{k-1}^{AGH} =D_{k-1}^A(I_{{}_N}+D_{k-1}^GD_{k-1}^H)^{-1}, &&
 D_{k-1}^{A^{\top}HG}= (D_{k-1}^A)^{\top}(I_{{}_N}+D_{k-1}^HD_{k-1}^G)^{-1}.
\end{array}
\right.
\] %\end{equation}

The corresponding low-rank factors are
\begin{eqnarray}
&&
\begin{array}{c@{\hspace{-2pt}}l}
\begin{array}{rr}
 \hspace{1.2cm} m^g_{k-1} \ \ \ \  m^{a_1}_{k-1} \hspace{1.1cm} m^{g}_{k-1} \hspace{1.3cm} \ \ m^{h}_{k-1}   \hspace{1.5cm} \ m^{a_2}_{k-1} \ \ \  &
\end{array}
& \vspace{0.1cm} \\
L_{k}^G = \left[
\begin{array}{ccccc}
L_{k-1}^G, & L_{1,{k-1}}^A, & D_{k-1}^{AGH}L_{k-1}^G, & D_{k-1}^{AGHG}L_{{k-1}}^H, & D_{k-1}^{AGHG}L_{2,{k-1}}^A
\end{array}
\right]
&
\begin{array}{l}
N,
\end{array}
\end{array}  \label{lgk}\\
&&
\begin{array}{c@{\hspace{-2pt}}l}
\begin{array}{rr}
 \hspace{2.4cm}   m^{a_1}_{k-1} \hspace{1.1cm} m^{g}_{k-1} \hspace{1.3cm} \ \ m^{h}_{k-1}   \hspace{1.5cm} \ m^{a_1}_{k-1} \ \ \  &
\end{array}
& \vspace{0.1cm} \\
L_{1,k}^A = \left[
\begin{array}{ccccc}
\ \ \ \ \ \ \ &L_{1,{k-1}}^A, & D_{k-1}^{AGH}L_{{k-1}}^G, &  D_{k-1}^{AGHG}L_{k-1}^H, & D_{k-1}^{AGH}L_{1,{k-1}}^A
\end{array}
\right]
&
\begin{array}{l}
N,
\end{array}
\end{array}\label{la1k}\\
&&
\begin{array}{c@{\hspace{-2pt}}l}
\begin{array}{rr}
 \hspace{.5cm} m^h_{k-1} \ \ \ \  m^{a_2}_{k-1} \hspace{1.1cm} m^{h}_{k-1} \hspace{1.3cm} \ \ m^{g}_{k-1}   \hspace{1.5cm} \ m^{a_1}_{k-1} \ \ \  &
\end{array}
& \vspace{0.1cm} \\
L_{k}^H = \left[
\begin{array}{ccccc}
L_{k-1}^H, & L_{2,{k-1}}^A, & D_{k-1}^{A^\top HG}L_{k-1}^H, & D_{k-1}^{A^\top HGH}L_{{k-1}}^G, & D_{k-1}^{A^\top HGH}L_{1,{k-1}}^A
\end{array}
\right]
&
\begin{array}{l}
N,
\end{array}
\end{array} \label{lhk}\\
&&
\begin{array}{c@{\hspace{-2pt}}l}
\begin{array}{rr}
 \hspace{1.8cm}   m^{a_2}_{k-1} \hspace{1.1cm} m^{g}_{k-1} \hspace{1.3cm} \ \ m^{h}_{k-1}   \hspace{1.5cm} \ m^{a_2}_{k-1} \ \ \  &
\end{array}
& \vspace{0.1cm} \\
L_{2,k}^A = \left[
\begin{array}{ccccc}
\ \ \ \ \ \ \ \ & L_{2,{k-1}}^A, & D_{k-1}^{A^\top HG}L_{{k-1}}^H, & D_{k-1}^{A^\top HGH}L_{k-1}^G, & D_{k-1}^{A^\top HG}L_{2,{k-1}}^A
\end{array}
\right]
&
\begin{array}{l}
N.
\end{array}
\end{array} \label{la2k}
\end{eqnarray}

To express kernels explicitly, let
\begin{eqnarray}
%\left.
%\begin{array}{rcl}
\Theta_{k-1}^H = (L_{{k-1}}^H)^\top D_{k-1}^{GHG}L_{{k-1}}^H,
& \Theta_{k-1}^G = (L_{{k-1}}^G)^\top D_{k-1}^{HGH}L_{{k-1}}^G, &
 \Theta_{k-1}^{HG} = (L_{{k-1}}^H)^\top D_{k-1}^{GH}L_{{k-1}}^G,  \nonumber\\%\label{tht-hgk}\\
\Theta_{{k-1}}^{A} = (L_{2,{k-1}}^A)^\top D_{k-1}^{GH}L_{1,{k-1}}^A,
&  \Theta_{1,{k-1}}^{A} = (L_{1,{k-1}}^A)^\top D_{k-1}^{HGH}L_{1,{k-1}}^A, &
\Theta_{2,{k-1}}^{A} = (L_{2,{k-1}}^A)^\top D_{k-1}^{GHG}L_{2,{k-1}}^A \nonumber % \label{tht-aak}
% \end{array}
%\right.
\end{eqnarray}
and
\begin{eqnarray}
\Theta_{1,{k-1}}^{AH} = (L_{1,{k-1}}^A)^\top D_{k-1}^{HG}L_{{k-1}}^H, &&
 \Theta_{1,{k-1}}^{AG} = (L_{1,{k-1}}^A)^\top D_{k-1}^{HGH}L_{{k-1}}^G, \nonumber \\
\Theta_{2,{k-1}}^{AH} = (L_{2,{k-1}}^A)^\top D_{k-1}^{GHG}L_{{k-1}}^H, &&
 \Theta_{2,{k-1}}^{AG} = (L_{2,{k-1}}^A)^\top D_{k-1}^{GH}L_{{k-1}}^G \nonumber%\label{tht-ahagk}
 \end{eqnarray}
with
\[ %\begin{equation}\label{dghhgk}
\left.
\begin{array}{rcl}
D_{k-1}^{GHG} =(I_{{}_N}+D_{k-1}^GD_{k-1}^H)^{-1}D_{k-1}^G, &&
D_{k-1}^{HGH} = (I_{{}_N}+D_{k-1}^HD_{k-1}^G)^{-1}D_{k-1}^H,\\
D_{k-1}^{GH} = (I_{{}_N}+D_{k-1}^GD_{k-1}^H)^{-1},  &&
D_{k-1}^{HG} = (I_{{}_N}+D_{k-1}^HD_{k-1}^G)^{-1}.
\end{array}
\right.
\] %\end{equation}

Define kernel components
\begin{equation}\label{kghk}
\left.
\begin{array}{rcl}
K_{k-1}^{GH} = \left[
\begin{matrix}
0 & K_{k-1}^G  \\
K_{k-1}^H & 0
\end{matrix} \right] \
\Big(I_{{}_{m_{{k-1}}^h+m_{{k-1}}^g}} + \left[
\begin{matrix}
-\Theta_{k-1}^H & \Theta_{k-1}^{HG} \\
(\Theta_{k-1}^{HG})^\top & \Theta_{k-1}^{G}
\end{matrix} \right] \
\left[
\begin{matrix}
-K_{k-1}^H & 0 \\
0 & K_{k-1}^{G}
\end{matrix} \right]
\Big)^{-1},
\end{array}\right.
\end{equation}

\begin{equation}
\left.
\begin{array}{rcl}
K_{k-1}^{GHG} = K_{k-1}^{GH} \left[
\begin{matrix}
0 & I_{m_{{k-1}}^h}  \\
-I_{m_{{k-1}}^g} & 0
\end{matrix} \right],  \ \
K_{k-1}^{HGH} = \left[
\begin{matrix}
0 & -I_{m_{{k-1}}^h}  \\
I_{m_{{k-1}}^g} & 0
\end{matrix} \right] K_{k-1}^{GH}
\end{array}\right.
\end{equation}
and
\begin{equation}\label{kghgk}
\left.
\begin{array}{rcl}
K_{k-1}^{AGHG} &=& -K_{k-1}^{A}[
\Theta_{2,{k-1}}^{AG}, \ \Theta_{2,{k-1}}^{AH}] K_{k-1}^{GHG},  \\
K_{k-1}^{A^\top HGH} &=& -(K_{k-1}^{A})^\top[
\Theta_{1,{k-1}}^{AH}, \ \Theta_{1,{k-1}}^{AG}] K_{k-1}^{HGH},  \\
K_{k-1}^{AGHGA^\top} &=& K_{k-1}^{A}\Theta_{2,{k-1}}^{A}(K_{k-1}^{A})^\top+ K_{k-1}^{AGHG} [\Theta_{2,{k-1}}^{AG}, \Theta_{2,{k-1}}^{AH}]^\top (K_{k-1}^{A})^\top,\\
K_{k-1}^{A^\top HGHA} &=& (K_{k-1}^{A})^\top \Theta_{1,{k-1}}^{A} K_{k-1}^{A} + K_{k-1}^{A^\top HGH}[
\Theta_{1,{k-1}}^{AH}, \ \Theta_{1,{k-1}}^{AG}]^\top K_{k-1}^{A},\\
K_{k-1}^{AGH} &=& -K_{k-1}^{A}[
\Theta_{2,{k-1}}^{AG}, \ \Theta_{2,{k-1}}^{AH}] K_{k-1}^{GH},\\
K_{k-1}^{A^\top GH} &=& -(K_{k-1}^{A})^\top[
\Theta_{1,{k-1}}^{AH}, \ \Theta_{1,{k-1}}^{AG}] (K_{k-1}^{GH})^\top,\\
K_{k-1}^{AGHA} &=& K_{k-1}^{A} \Theta_{{k-1}}^{A} K_{k-1}^{A} + K_{k-1}^{AGH}[
\Theta_{1,{k-1}}^{AH}, \ \Theta_{1,{k-1}}^{AG}]^\top K_{k-1}^{A}.
\end{array}\right.
\end{equation}
Then kernel matrices corresponding to $L^G_k$,  $L^H_k$ and $L^A_{1,k}$ ($L^A_{2,k}$) at $k$-th step could be represented as
\begin{equation}\label{kgkp1}
\begin{array}{c@{\hspace{-2pt}}l}
\begin{array}{rr}
 \hspace{1.8cm} m^g_{k-1}  \ \ \ \ \ \ \ m^{a_1}_{k-1}   \ \ \  m^{g}_{k-1}+m^{h}_{k-1} \  m^{a_2}_{k-1}  &
\end{array}
& \vspace{0.1cm} \\
K_{k}^{G} = \left[
\begin{array}{cccc}
  K_{k-1}^G  & 0  &0 &0 \\
0& K_{k-1}^{AGHGA^\top} & K_{k-1}^{AGHG} & K_{k-1}^{A} \\
0 & (K_{k-1}^{AGHG})^\top & -K_{k-1}^{GHG} &0 \\
0 & (K_{k-1}^{A})^\top & 0 & 0
\end{array}
\right]
&
\begin{array}{l}
\hspace{-.3cm}m^g_{k-1} \\ \hspace{-.3cm}m^{a_1}_{k-1} \\ \hspace{-.3cm}m^{g}_{k-1}+m^{h}_{k-1} \\ \hspace{-.3cm}m^{a_2}_{k-1} \\
\end{array}
\end{array},
\end{equation}

\begin{equation}%\label{khkp1}
\begin{array}{c@{\hspace{-2pt}}l}
\begin{array}{rr}
 \hspace{1.3cm} m^h_{k-1}  \ \ \ \ \ \ \ m^{a_2}_{k-1}  \ \ \ \ \  m^{h}_{k-1}+m^{g}_{k-1} \ \ \ m^{a_1}_{k-1}  &
\end{array}
& \vspace{0.1cm} \\
K_{k}^{H} = \left[
\begin{array}{cccc}
  K_{k-1}^H  & 0  &0 &0 \\
0& K_{k-1}^{A^\top HGHA} & K_{k-1}^{A^\top HGH} & (K_{k-1}^{A})^\top \\
0 & (K_{k-1}^{A^\top HGH})^\top & -K_{k-1}^{HGH} &0 \\
0 & K_{k-1}^{A} & 0 & 0
\end{array}
\right]
&
\begin{array}{l}
m^h_{k-1} \\ m^{a_2}_{k-1} \\ m^{h}_{k-1}+m^{g}_{k-1} \\ m^{a_1}_{k-1} \\
\end{array}
\end{array}
\end{equation}
and

\begin{equation}\label{kakp1}
\begin{array}{c@{\hspace{-2pt}}l}
\begin{array}{rr}
 \hspace{2cm} m^{a_2}_{k-1} \ \ m^{g}_{k-1}+m^{h}_{k-1} \ m^{a_2}_{k-1} &
\end{array}
& \vspace{0.1cm} \\
K_{k}^{A} = \left[
\begin{array}{ccc}
K_{k-1}^{AGHA} & K_{k-1}^{AGH} &K_{k-1}^{A}   \\
(K_{k-1}^{A^\top GH})^\top & -K_{k-1}^{GH} & 0 \\
K_{k-1}^A  & 0 &0
\end{array}
\right]
&
\begin{array}{l}
m^{a_1}_{k-1} \\ m^{g}_{k-1}+m^{h}_{k-1} \\ m^{a_1}_{k-1}
\end{array}
\end{array}.
\end{equation}

%\begin{equation}\label{laghkpk}
%\left.
%\begin{array}{rcl}
%L_{1,k}^A &=& [ \ \ \ \ \ \ \ \ \  L_{1,{k-1}}^A, \ D_{k-1}^{AGH}L_{{k-1}}^G, \  D_{k-1}^{AGHG}L_{k-1}^H, \ D_{k-1}^{AGH}L_{1,{k-1}}^A \  ], \\
%L_{k}^H &=&[L_{k-1}^H, \ L_{2,{k-1}}^A, \ D_{k-1}^{A^\top HG}L_{k-1}^H, \ D_{k-1}^{A^\top HGH}L_{{k-1}}^G, \ D_{k-1}^{A^\top HGH}L_{1,{k-1}}^A], \\
%L_{2,k}^A &=&[ \ \ \ \ \ \ \ \ \ L_{2,{k-1}}^A, \ D_{k-1}^{A^\top HG}L_{{k-1}}^H, \ D_{k-1}^{A^\top HGH}L_{k-1}^G, \ D_{k-1}^{A^\top HG}L_{2,{k-1}}^A \ \ ].
%\end{array}
%\right.
%\end{equation}

\begin{remark}
1. The banded part \eqref{sdak} in FSDA can iterate independently of the low-rank part,  contributing to the motivation to establish the pre-termination condition in Section 5.

2. Low-rank factors in \eqref{lgk}-\eqref{la2k} are seen growing at least with a scale of $O(4^k)$, obviously intolerable for large-scale problems. So a deflation process and a truncation and compression technique are required to reduce the column dimension of low-rank factors as much as possible.

3. In real implementations, low-rank factors and kernels for $k\geq 2$ are actually deflated, truncated and compressed ones as described in the next two sections, where a superscript ``$dt$'' is labelled in the upper right corner of each low-rank factor. Correspondingly, column numbers $m^g_{k-1}$, $m^h_{k-1}$, $m^{a_1}_{k-1}$ and $m^{a_2}_{k-1}$ are also ones after deflation, truncation and compression. Here we temporarily omit this superscript ``$dt$'' just for the convenience of describing the successive iteration process.

\end{remark}

\subsection{Convergence of factors and kernels.}
The convergence of the sequences $\{G_k\}$, $\{H_k\}$ and $\{A_k\}$ in SDA \eqref{sda} is given in \cite[Thm~3.1]{lx06}.

\begin{theorem} \label{thm21}
Assume that $X,Y >0$ satisfy the DARE \eqref{d} and its dual equation \eqref{dd} and let
\[
S := (I+GX)^{-1} A, \ \ T:= (I+HY)^{-1} A^\top.
\]
Then the matrix sequences $\{A_k\}$, $\{G_k\}$ and $\{H_k\}$ generated by the SDA satisfy
\begin{itemize}
 \item[(1)] $A_k = (I + G_k X)S^{2^k}$;
\item[(2)] $H \leq H_k \leq H_{k+1} \leq X$ and
$X - H_k = \left( S^\top \right)^{2^k} (X + XG_k X)S^{2^k} \leq \left( S^\top \right)^{2^k} (X + XY X)S^{2^k}$;
\item[(3)] $G \leq G_k \leq G_{k+1} \leq Y$ and
$Y - G_k = \left( T^\top \right)^{2^k} (Y + YH_k Y)T^{2^k} \leq \left( T^\top \right)^{2^k} (Y + YX Y)T^{2^k}$.
\end{itemize}
\end{theorem}

Let $X = D^X+ L^XK^X(L^X)^\top$ and $Y=D^Y+ L^YK^Y(L^Y)^\top$ be the stabilizing solutions of DARE \eqref{d} and its dual equation \eqref{dd}, respectively. It follows from (1) of the above theorem that
$\|A_k\| \leq  (1+\|X\|\cdot \|Y\|)\|S^{2^k}\|$, indicating that $A_k$ converges to zero quadratically when $\rho(S)<1$.
By recalling the decomposition $A_k = D^A_k + L^A_{1,k}K^A_k(L^A_{2,k})^\top$, the banded matrix sequence $\{D^A_k\}$ and the low-rank sequence $\{L^A_{1,k}K^A_kL^A_{2,k}\}$ will converge to zero respectively. Similarly, it follows form (2) and (3) that
$\|H_k -X\| \leq  \|X\|(1+\|X\|\cdot \|Y\|)\| S^{2^k}\|^2$ and $\|G_k - Y\| \leq  \|Y\|(1+\|X\|\cdot \|Y\|)\| T^{2^k}\|^2$, indicating that $H_k$ and $G_k$ respectively converge to $X$ and $Y$ quadratically when $\rho(S)<1$ and $\rho(T)<1$. Then from the iterative decomposition $H_k = D^H_k+L^H_k K^H_k (L^H_k)^\top$ and $G_k = D^G_k+L^G_k K^G_k (L^G_k)^\top$, the banded sequences $\{D^H_k\}$ and $\{D^G_k\}$ will respectively converge to banded parts $D^X$ and $D^Y$ of $X$ and $Y$. Also, the low-rank sequences $\{L^H_k\}$, $\{L^G_k\}$ and the kernel sequences $\{K^H_k\}$, $\{K^G_k\}$ will converge to low-rank parts $L^X$, $L^Y$ and kernels $K^X$, $K^Y$ of $X$ and $Y$, respectively. We can conclude the above as the following corollary.

\begin{corollary}%\label{cor21}
Suppose that $X = D^X+ L^XK^X(L^X)^\top$ and $Y=D^Y+ L^YK^Y(L^Y)^\top$ satisfy the DARE and its dual equation, respectively. If $\rho(S)<1$ and $\rho(T)<1$, then for FSDA, the sequences $\{D^A_k\}$ and  $\{L^A_{1,k}K^A_kL^A_{2,k}\}$ converge to zero with
\[
\limsup_{k\rightarrow\infty}\sqrt[2^k]{\|D^A_k\|} \leq \rho(S),  \ \ \limsup_{k\rightarrow\infty}\sqrt[2^k]{\|L^A_{1,k}K^A_kL^A_{2,k}\|}\leq \rho(T).
\]
The banded sequences $\{L^H_k\}$, $\{L^G_k\}$, low-rank sequences $\{L^H_k\}$, $\{L^G_k\}$ and the kernel sequences $\{K^H_k\}$, $\{K^G_k\}$ will respectively converge to $D^X$, $D^Y$, $L^X$, $L^Y$ and the kernel $K^X$, $K^Y$ with
\begin{eqnarray}
\limsup_{k\rightarrow\infty}\sqrt[2^k]
{\|D^H_k-D^X\|}\leq
\rho^2(S), &&\limsup_{k\rightarrow\infty}\sqrt[2^k] {\|D^G_k-D^Y\|}\leq
\rho^2(T), \nonumber\\
\limsup_{k\rightarrow\infty}\sqrt[2^k]
{\|L^H_k-L^X\|}\leq
\rho^2(S), && \limsup_{k\rightarrow\infty}\sqrt[2^k]{\|L^G_k-L^Y\|}\leq
\rho^2(T), \nonumber \\
\limsup_{k\rightarrow\infty}\sqrt[2^k]{\|K^H_k-K^X\|}\leq
\rho^2(S), && \limsup_{k\rightarrow\infty}\sqrt[2^k]{\|K^G_k-K^Y\|}\leq
\rho^2(T). \nonumber
\end{eqnarray}
\end{corollary}

\begin{remark}
 Although the product $L^A_{1,k}K^A_{k}L^A_{2,k}$ converges to zero quadratically, it follows from \eqref{la1k}, \eqref{la2k} and \eqref{kakp1} that the kernel $K^A_k$ and low-rank factors $L^A_{1,k}$ and $L^A_{2,k}$ might still not converge to zero, respectively.
\end{remark}

\section{Deflation of low-rank factors and kernels}
It has been shown that there is an exponential increase in columns of low-rank factors and kernels. Nevertheless, it is clear that the first three items in $L_{1,k}^A$ and $L_{2,k}^A$ (see \eqref{la1k}, \eqref{la2k}) are same to the second to the fourth item in $L_{k}^G$ and $L_{k}^H$(see \eqref{lgk}, \eqref{lhk}), respectively. Moreover, a careful observation on factor matrices reveals that there are some items in $L^G_k$ and $L^H_k$ are essentially repetitive to some in $L^A_{1,k}$ and $L^A_{2,k}
$, respectively. Then the deflation of low-rank factors and kernels are definitely required. To see this process clearly, we start with the case $k=2$.

\subsubsection*{Case for $k=2$.}
Consider the deflation of low-rank factors firstly. It follows from \eqref{lgk}-\eqref{la2k} that
\[%\begin{equation}\label{lagh2}
\left.
\begin{array}{rcl}
L_2^G &=&[L_1^G, \ L_{11}^A, \ D_1^{AGH}L_1^G, \ D_1^{AGHG}L_{1}^H, \ D_1^{AGHG}L_{21}^A], \\
L_{12}^A &=& [L_{11}^A, \ D_1^{AGH}L_{1}^G, \  D_1^{AGHG}L_1^H, \ D_1^{AGH}L_{11}^A  ], \\
L_2^H &=&[L_1^H, \ L_{21}^A, \ D_1^{A^\top HG}L_1^H, \ D_1^{A^\top HGH}L_{1}^G, \ D_1^{A^\top HGH}L_{11}^A], \\
L_{22}^A &=& [ L_{21}^A, \ D_1^{A^\top HG}L_{1}^H, \ D_1^{A^\top HGH}L_1^G, \ D_1^{A^\top HG}L_{21}^A ]
\end{array}
\right.
\]%\end{equation}
with
\[%\begin{equation}\label{daghahg1}
\left.
\begin{array}{rcl}
D_1^{AGHG} =D_1^A(I+D_1^GD_1^H)^{-1}D_1^G, && D_1^{A^{\top}HGH}= (D_1^A)^{\top}D_1^H(I+D_1^GD_1^H)^{-1},\\
 D_1^{AGH} =D_1^A(I+D_1^GD_1^H)^{-1}, &&
 D_1^{A^{\top}HG}= (D_1^A)^{\top}(I+D_1^HD_1^G)^{-1}.
\end{array}
\right.
\]%\end{equation}
Expanding the above low-rank factors with initial $L^A_{10}\in {\mathbb R}^{N\times m^a}$ and $L^A_{20}\in {\mathbb R}^{N\times m^a}$, one can see from Appendix A that $L_{10}^A$ and $D_1^{AGHG}L_{20}^A$ (or $L_{20}^A$ and $D_1^{A^\top HGH}L_{10}^A$) occur twice in
$L_2^{G}$ (or $L_2^{H}$). To reduce the column dimension of $L^G_k$, we will
shift the overlapped $L^A_{10}$ in $L^G_1$ to the one in $L^A_{11}$ and $D_1^{AGHG}L_{20}^A$ in $D_1^{AGHG}L_{21}^A$ to the one in $D_1^{AGHG}L_{1}^H$, respectively. Then the original $L^G_2$ is deflated to $L^{Gd}_2$ of the smaller column dimension, where the superscript ``{\it d} '' means the matrix after deflation. Analogously, as
$D_1^{AGH}L_{10}^A$ (or $D_1^{A^\top HG}L_{20}^A$) appears twice in $L_{12}^{A}$ (or $L_{22}^{A}$), the overlapped $L^A_{20}$ in $L^H_1$ will be shifted to the one in $L^A_{21}$ and $D_1^{A^\top HGH}L_{10}^A$ in $D_1^{A^\top HGH}L_{11}^A$ to the one in $D_1^{A^\top HGH}L_{1}^G$, respectively. Then the original $L^H_2$ is deflated to $L^{Hd}_2$. Such a process can be applied to $L^A_{12}$ and $L^A_{22}$, obtaining $L^{Ad}_{12}$ and $L^{Ad}_{22}$ listed in Appendix A, respectively, where the left blank in each factor corresponds the deleted matrix and the black bold  matrices inhere from the un-deflated ones.

For kernels at $k=2$, one has
\[%\begin{equation}\label{kg2}
\begin{array}{c@{\hspace{-2pt}}l}
\begin{array}{rr}
 \hspace{1.3cm} 2m^a \ \  \ \ \ \ \ \ 4m^{a} \ \ \ \ \ \ \ \ \ \ 2m^{a} \ \  \ \ \ \ 2m^{a}   &
\end{array}
& \vspace{0.1cm} \\
K_{2}^{G} = \left[
\begin{array}{cccc}
  K_{1}^G  & 0  &0 &0 \\
0& K_{1}^{AGHGA^\top} & K_{1}^{AGHG} & K_{1}^{A} \\
0 & (K_{1}^{AGHG})^\top & -K_{1}^{GHG} &0 \\
0 & (K_{1}^{A})^\top & 0 & 0
\end{array}
\right]
&
\begin{array}{l}
2m^a \\ 4m^{a} \\ 2m^{a} \\ 2m^{a} \\
\end{array}
\end{array},
\]%\end{equation}

\[%\begin{equation}\label{kh2}
\begin{array}{c@{\hspace{-2pt}}l}
\begin{array}{rr}
 \hspace{1.3cm} 2m^a \ \ \ \ \ \ \ \ \ 4m^{a} \ \ \ \ \ \ \ \ \ \ \ \ 2m^{a} \ \ \ \ \ \ \ \ 2m^{a}   &
\end{array}
& \vspace{0.1cm} \\
K_{2}^{H} = \left[
\begin{array}{cccc}
  K_{1}^H  & 0  &0 &0 \\
0& K_{1}^{A^\top HGHA} & K_{1}^{A^\top HGH} & (K_{1}^{A})^\top \\
0 & (K_{1}^{A^\top HGH})^\top & -K_{1}^{HGH} &0 \\
0 & K_{1}^{A} & 0 & 0
\end{array}
\right]
&
\begin{array}{l}
2m^a \\ 4m^{a} \\ 2m^{a} \\ 2m^{a} \\
\end{array}
\end{array}
\]%\end{equation}
and
\[%\begin{equation}\label{ka2}
\begin{array}{c@{\hspace{-2pt}}l}
\begin{array}{rr}
 \hspace{2cm} 2m^{a} \ \ \ \ \ \ \ \ \  4m^{a}\ \ \ \ \ 2m^{a} &
\end{array}
& \vspace{0.1cm} \\
K_{2}^{A} = \left[
\begin{array}{ccc}
K_{1}^{AGHA} & K_{1}^{AGH} &K_{1}^{A}   \\
(K_{1}^{A^\top GH})^\top & -K_{1}^{GH} & 0 \\
K_{1}^A  & 0 &0
\end{array}
\right]
&
\begin{array}{l}
2m^{a} \\ 4m^{a} \\ 2m^{a}
\end{array}
\end{array}
\]%\end{equation}
with non-zero components defined in \eqref{kghk}-\eqref{kghgk}. Here details of the deflation of $K_2^G$ is explained explicitly and the implementation of $K_2^H$ is the same. In fact, there are 10 block rows and block columns with each of initial size $m^a\times m^a$ in $K_2^G$. The
deflation process will simultaneously shift the first block row and block column to the third block row and block column and the ninth block row and block column to the seventh block row and block
column. Then the (2,2) sub-block of $K_1^G$ covers the (1,1) sub-block of $K_1^{AGHGA^\top}$. The first column sub-block of $K_1^{A}$ and the first row sub-block of $(K_1^{A})^\top$ overlap with the third column sub-block of $K_1^{AGHA}$ and the third row sub-block of $(K_1^{AGHA})^\top$, respectively, completing the deflated matrix $K^{Gd}_2$.

\begin{figure}[H]
\centering
$K_{2}^G$:
\subfigure{
\begin{minipage}{4.6cm}
\includegraphics[height=4.5cm, width =4.7cm]{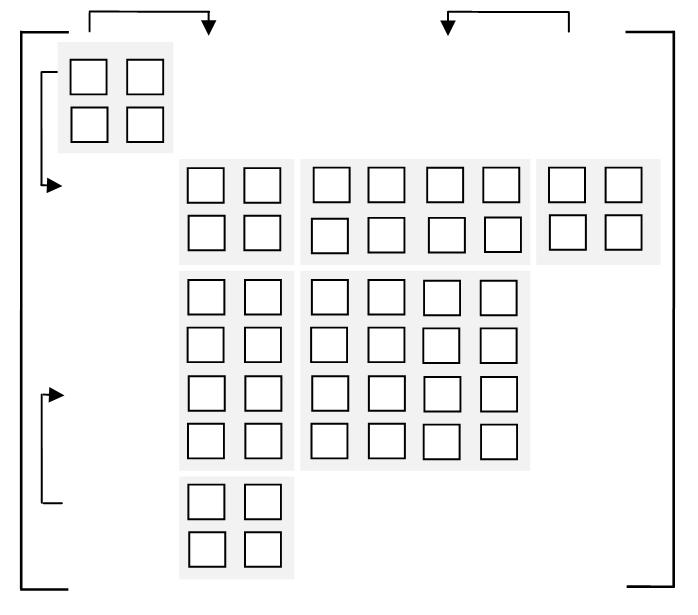}
\end{minipage}
}
$\stackrel{d}{\rightarrow}$
\begin{minipage}{3.2cm} %[b]%{0.2\textwidth}
\subfigure{
\includegraphics[height=3.0cm, width =3.2cm]{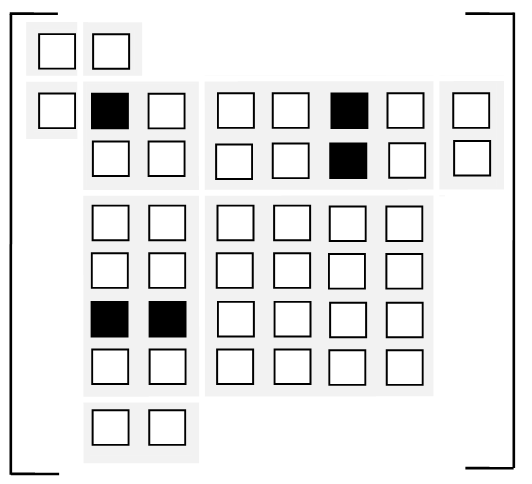}
}
\end{minipage}
$:=K^{Gd}_2$
\begin{center}
Fig 3.1. The deflation process of $K_2^G$ (or $K_2^H$).
\end{center}
\end{figure}

\begin{figure}[H]
\centering
$K_{2}^A$:
\subfigure{
\begin{minipage}{4.2cm}
\includegraphics[height=4.0cm, width =4.2cm]{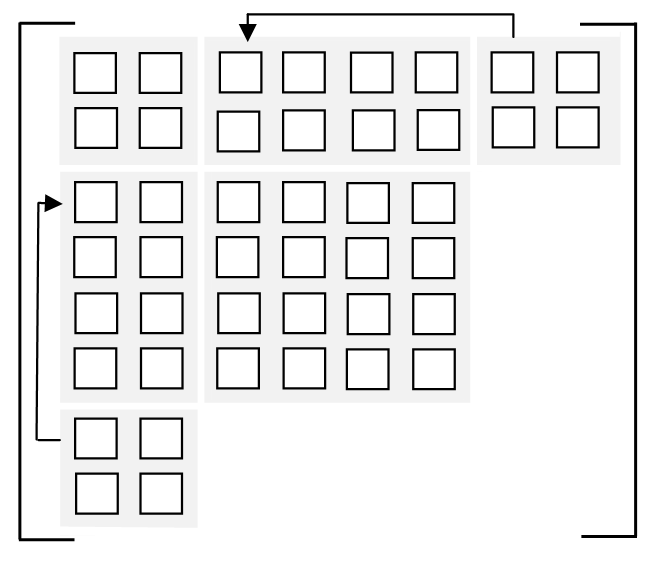}
\end{minipage}
}
$\stackrel{d}{\rightarrow}$
\begin{minipage}{3.1cm} %[b]%{0.2\textwidth}
\subfigure{
\includegraphics[height=2.9cm, width =3.1cm]{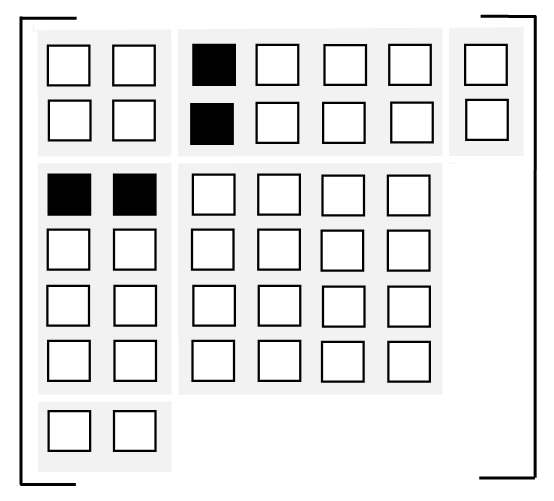}
}
\end{minipage}
$:=K^{Ad}_2$
\begin{center}
Fig 3.2. The deflation process of $K_2^A$.
\end{center}
\end{figure}

Analogously, there
are 8 block rows and block columns with each of the initial size $m^a\times m^a$ in $K_2^A$. The deflation process simultaneously
shifts the seventh column and row sub-blocks to the third column and row sub-blocks, respectively. Then the first column sub-block of the upper right $K^A_1$ and the first row sub-block of the lower-left $K_1^{A}$ overlap with the first column sub-block of $K_1^{AGH}$ and the first row sub-block of $(K_1^{A^\top GH})^\top$, respectively, completing the deflation of $K^{Ad}_2$.

The whole process is described in Fig. 3.1 and 3.2 where each small square is of size $m^a\times m^a$ and the  grey-bottom block represents the non-zero component in $K_2^G$ and $K_2^A$. The little white squares in $K^{Gd}_2$ and $K^{Ad}_2$ inhere from the originally un-deflated sub-matrices and the little black squares in $K^{Gd}_2$ and $K^{Ad}_2$ represent the overlayed sub-matrices.

\subsubsection*{Case for $k\geq 3$.}
After the $(k-1)$-th deflation, there are $m^g_{k-1}-km^a$ (or $m^h_{k-1}-km^a$) columns in $L^G_{k-1}$ and $L^A_{1,k-1}$ (or $L^H_{k-1}$ and $L^A_{2,k-1}$) and
$m^{a_2}_{k-1}-m^a$ (or $m^{a_1}_{k-1}-m^a$) columns in $D^{AGHG}_{k-1}L^A_{2,k-1}$ and $D^{AGHG}_{k-1}L^H_{k-1}$ (or $D^{A^\top HGH}_{k-1}L^A_{1,k-1}$ and $D^{A^\top HGH}_{k-1}L^G_{k-1}$) are identical. Then, one can shift columns of
$$L^G_{k-1}(:, (k-1)m^a+1:
m^g_{k-1}-m^a) \ \ \Big(\mbox{or} \ L^H_{k-1}(:, (k-1)m^a+1: m^h_{k-1}-m^a)\ \Big)$$
to columns of
$$L^A_{1,k-1}(:, 1: m^g_{k-1}-km^a) \ \ \Big(\mbox{or} \ L^A_{2,k-1}(:, 1: m^h_{k-1}-km^a)\ \Big)$$
and columns of
$$D^{AGHG}_{k-1}L^A_{2,k-1}(:, 1:m^{a_2}_{k-1}-m^a)\ \ \Big(\mbox{or} \ D^{A^\top HGH}_{k-1}L^A_{1,k-1}(:, 1:m^{a_1}_{k-1}-m^a)\ \Big)$$
to columns of
$$D^{AGHG}_{k-1}L^H_{k-1}(:, m^h_{k-1}-m^{a_2}_{k-1}+1: m^h_{k-1}-m^a) \ \ \Big(\mbox{or} \ D^{A^\top HGH}_{k-1}L^G_{k-1}(:, m^g_{k-1}-m^{a_1}_{k-1}+1:
m^g_{k-1}-m^a)\ \Big)$$
in $L^G_k$ (or $L^H_k$), respectively.
Now, there are
$k-1$ matrices with each of order $N\times m^a$ are left in $L^G_{k-1}$ and $L^H_{k-1}$ (see the first item in \eqref{lgdk} and \eqref{lhdk} of Appendix B). Meanwhile, only one matrix of order $N\times m^a$ is left in $D^{AGHG}_kL^A_{2,k}
$, $D^{AGH}_kL^A_{1,k}$, $D^{A^\top HGH}_kL^A_{1,k}$ and $D^{A^\top HG}_kL^A_{2,k}
$(see the last item in \eqref{lgdk}--\eqref{lad2k} of Appendix B).

To deflate $L^A_{1,k}$ ($L^A_{2,k}$), columns of
$$D^{AGH}_{k-1}L^A_{1,k-1}(:, 1:m^{a_1}_{k-1}-m^a) \ \ \Big(\mbox{or} \ D^{A^\top HG}_{k-1}L^A_{2,k-1}(:,\ 1:m^{a_2}_{k-1}-m^a)\ \Big)$$ are shifted to columns of
$$D^{AGH}_{k-1}L^G_{k-1}(:, m^g_{k-1}-m^{a_1}_{k-1}+1: m^g_{k-1}-m^a)\ \ \Big(\mbox{or} \ D^{A^\top HG}_{k-1}L^H_{k-1}(:,\ m^h_{k-1}-m^{a_2}_{k-1}+1: m^h_{k-1}-m^a)\ \Big)$$ in $L^A_{1,k}$ (or $L^A_{2,k}$).
Appendix B gives details of the entire process, where low-rank factors in the $(k-1)$-th iteration are actually ones after deflation, truncation and compression, but only the superscript ``$d$'' is used just for the convenience of interpretation of the deflation process.

Correspondingly, kernel matrices $K^G_k$, $K^H_k$ and $K^A_k$ are deflated according to low-rank factors. Here we describe the deflation of $K_k^G$ and the way of $K_k^H$ is almost the same. By recalling the place of non-zero sub-matrices (the grey-bottom block in Fig. 3.3) of $K^G_k$ in \eqref{kgkp1}, the
deflation process essentially shifts $K^G_{k-1}((k-1)m^a+1: m^g_{k-1}-m^a, \ (k-1)m^a+1: m^g_{k-1}-m^a)$ to $K^{AGHGA^\top}_{k-1}(1:m^g_{k-1}-km^a, \ 1:m^g_{k-1}-km^a)$, columns $K^A_{k-1}(:, 1:m^{a_2}_{k-1}-m^a)$ to $K^{AGHG}_{k-1}(:,
m^{g}_{k-1}+m^{h}_{k-1}-m^{a_2}_{k-1}+1:m^{g}_{k-1}+m^h_{k-1}-m^a)$ and rows $(K^A_{k-1})^\top(1:m^{a_2}_{k-1}-m^a,\ :)$ to
$(K^{AGHG}_{k-1})^\top(m^{g}_{k-1}+m^{h}_{k-1}-m^{a_2}_{k-1}+1:m^{g}_{k-1}+m^h_{k-1}-m^a,\ :)$, respectively.See Figure 3.3 for more details.

\begin{figure}[H]
\centering
$K_{k}^G$:
\subfigure{
\begin{minipage}{6.3cm}
\includegraphics[height=6cm, width =6.3cm]{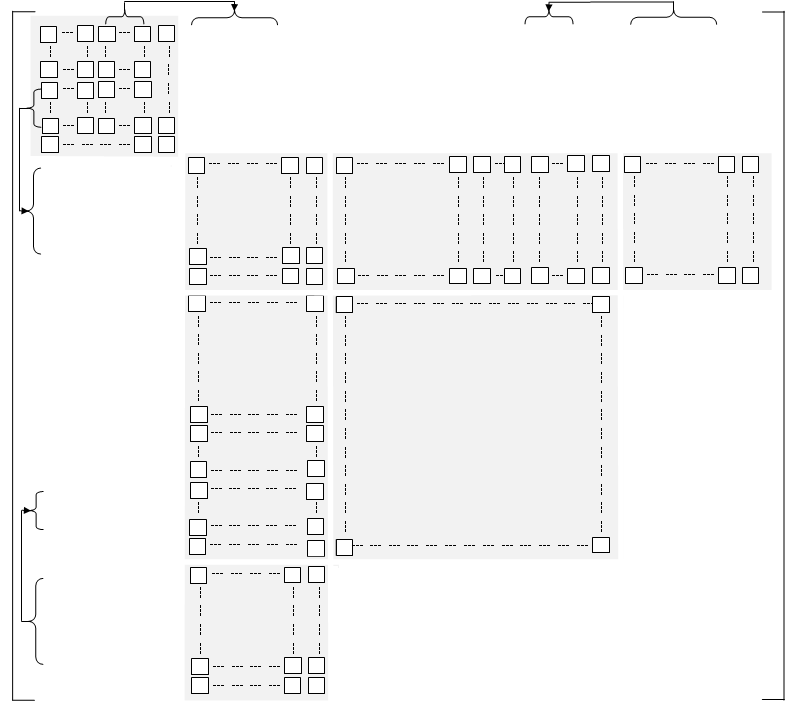}
\end{minipage}
}
$\stackrel{d}{\rightarrow}$
\begin{minipage}{5.2cm} %[b]%{0.2\textwidth}
\subfigure{
\includegraphics[height=5.2cm, width =5.2cm]{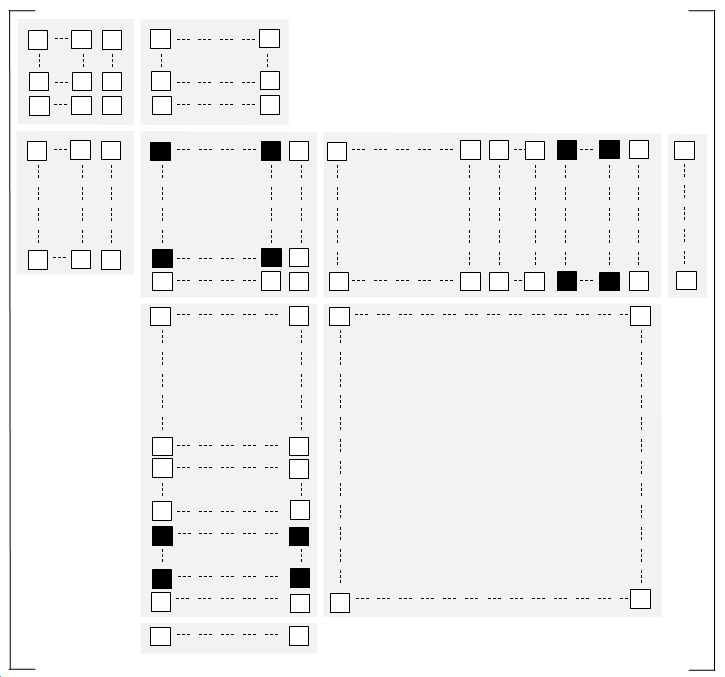}
}
\end{minipage}
$:=K_{k}^{Gd}$
\begin{center}
Fig. 3.3.  The deflation process of $K_{k}^G$ (or $K_{k}^H$).
\end{center}
\end{figure}

Similarly, by recalling the place of non-zero matrices (the grey-bottom block in Fig. 3.4) of $K^A_k$ in \eqref{kakp1}, the deflation process will shift columns $K^A_{k-1}(:, 1:m^{a_2}_{k-1}-m^a)$ to columns $K^{AHG}_{k-1}(:, m^h_{k-1}-m^{a_2}_{k-1}+1:m^h_{k-1}-m^a)$ and rows $K^A_{k-1}(1:m^{a_1}_{k-1}-m^a,:)$ to rows $(K^{A^\top GH}_{k-1})^\top( m^{g}_{k-1}-m^{a_1}_{k-1}+1:m^{g}_{k-1}-m^a,:)$.See Figure 3.4 for more details.

\begin{figure}[H]
\centering
$K_{k}^A$:
\subfigure{
\begin{minipage}{6.1cm}
\includegraphics[height=6.1cm, width =6.1cm]{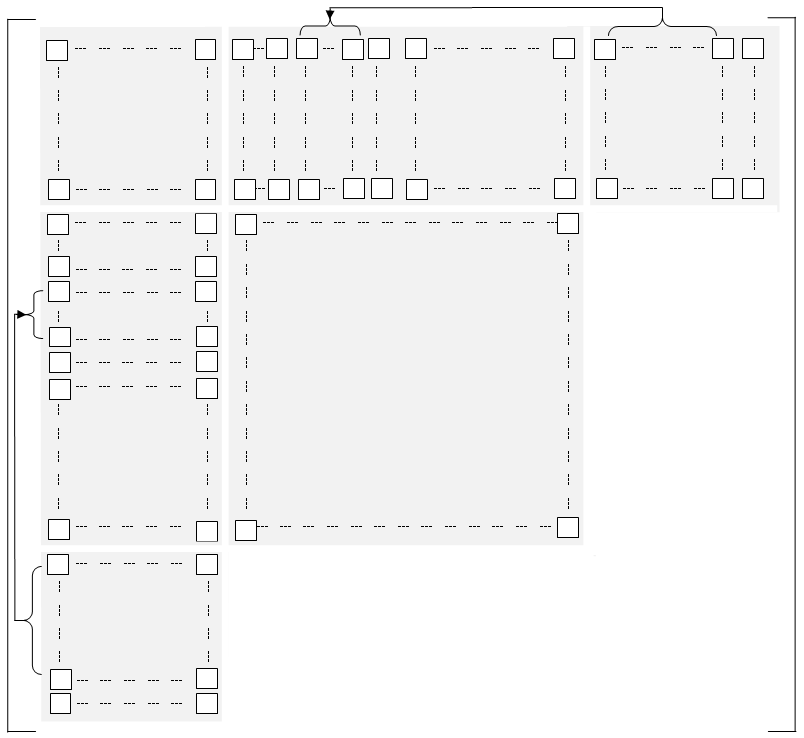}
\end{minipage}
}
$\stackrel{d}{\rightarrow}$
\begin{minipage}{5.2cm} %[b]%{0.2\textwidth}
\subfigure{
\includegraphics[height=5.2cm, width =5.2cm]{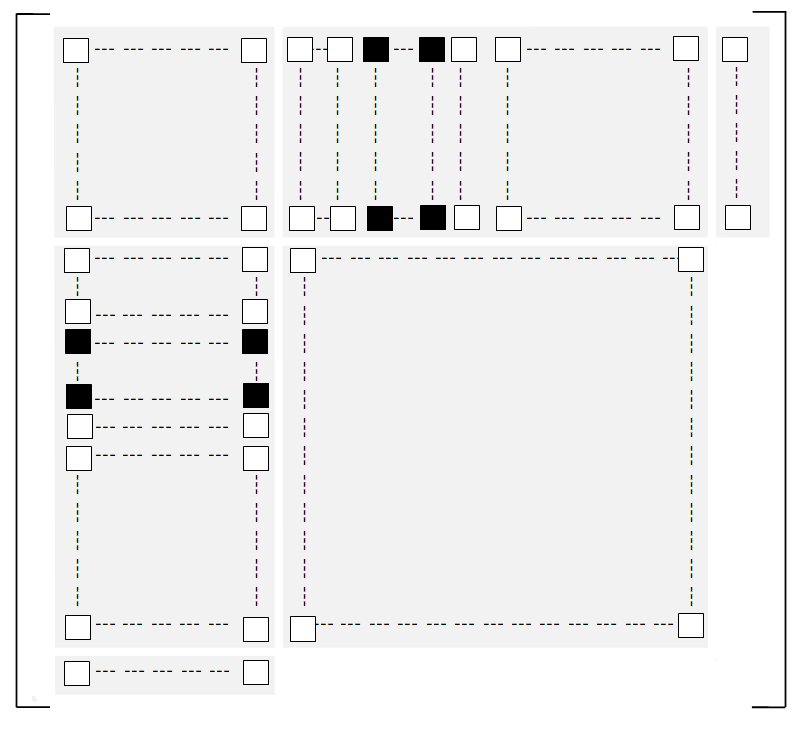}
}
\end{minipage}
$:=K_{k}^{Ad}$
\begin{center}
Fig. 3.4.  The deflation process of $K_{k}^A$.
\end{center}
\end{figure}

\section{Partial truncation and compression}
Although deflation of low-rank factors and kernels in last subsection can reduce dimensional growth of columns, the exponential increment of the un-deflated part is still rapid,  making large-scale computation and storage infeasible. Conventionally, one efficient way to shrink the column number of low-rank factors is the technique of truncation and compression (TC) in \cite{cw15, yfc19}, which unfortunately, is hard to be applied to our case due to the following two main obstacles.
\begin{itemize}
  \item Direct application of TC to $L^{Hd}_k$, $L^{Gd}_k$, $L^{Ad}_{1,k}$, $L^{Ad}_{2,k}$ and their corresponding kernels $K^{Hd}_k$, $K^{Gd}_k$ and $K^{Ad}_k$ at the $k$-th step will require four QR decompositions, resulting in a relatively high computational complexity and CPU consumption.

  \item The TC process applied to the whole low-rank factors at current step breaks up the implicit structure, causing  the deflation to be unrealized in the next iteration.
\end{itemize}

In this section, we will instead present a partial truncation and compression (PTC) to conquer the above two difficulties. Our PTC only requires two QR decompositions of the exponential increasing (not the entire) parts of low-rank factors and is capable of keeping the successive deflation for subsequent iterations.

\vspace{.3cm}
{\bf PTC for low-rank factors.} Recall the deflated forms \eqref{lgdk} and \eqref{lhdk} in Appendix B. $L^{Gd}_{k}$ and $L^{Hd}_{k}$ can be divided to three parts
 \[
\left.
\begin{array}{l}
L^{Gd}_{k} =[L^{Gd}_{k}(1), \ L^{Gd}_{k}(2), \ L^{Gd}_{k}(3)]\\
L^{Hd}_{k}=[L^{Hd}_{k}(1), \ L^{Hd}_{k}(2), \ L^{Hd}_{k}(3)].
\end{array}
\right.
\]
The first parts
\[
\left.
\begin{array}{l}
L^{Gd}_{k}(1):=[D_0^{AGHG}L_{20}^A, D_1^{AGHG}D_0^{A^\top GH}L_{20}^A,..., D_{k-2}^{AGHG}\Pi_{i=0}^{k-3}D_i^{A^\top GH}L_{20}^A]\in \mathbb{R}^{N \times (k-1)m^a}
\end{array}
\right.
\]
and
\[
\left.
\begin{array}{l}
L^{Hd}_{k}(1):=[D_0^{A^\top HGH}L_{10}^A, D_1^{A^\top HGH}D_0^{A^\top GH}L_{10}^A,..., D_{k-2}^{A^\top HGH}\Pi_{i=0}^{k-3}D_i^{AGH}L_{10}^A]\in \mathbb{R}^{N \times (k-1)m^a}
\end{array}
\right.
\]
rise only linearly with $k$, respectively, and the last parts
\[
\left.
\begin{array}{l}
L^{Gd}_{k}(3):= D_{k-1}^{AGHG}\Pi_{i=0}^{k-2}D_i^{A^\top GH}L_{20}^A\in \mathbb{R}^{N \times m^a}
\end{array}
\right.
\]
and
\[
\left.
\begin{array}{l}
L^{Hd}_{k}(3):= D_{k-1}^{A^\top HGH}\Pi_{i=0}^{k-2}D_i^{A GH}L_{10}^A\in \mathbb{R}^{N \times m^a}
\end{array}
\right.
\]
remain unchanged about the size $N\times m^a$, respectively. So we only truncate and compress the dominantly growing parts
\[
\left.
\begin{array}{l}
L^{Gd}_{k}(2):=[L^A_{1,k-1}, \ D^{AGH}_{k-1}L^G_{k-1}, \ D^{AGHG}_{k-1}L^H_{k-1}]
\end{array}
\right.
\]
and
\[
\left.
\begin{array}{l}
L^{Hd}_{k}(2):=[L^A_{2,k-1}, \ D^{A^\top HG}L^H_{k-1}, \ D^{A^\top HGH}_{k-1}L^G_{k-1}]
\end{array}
\right.
\]
by orthogonalization.
%We only describe the process for $L^H_k$ and $R^H_k$, operations on $L^G_k$, $L^A_{1,k}$, $L^A_{2,k}$, $R^G_k$ and $R^G_k$ are the same.
%Now we recall an important aspect of the SDA for large-scale DAREs --- the growth of $B_k$. As $H_k = M_k + C_k^\top T_k C_k^\top$ is not computed explicitly %\footnote{In the original SDA for small-scale AREs \cite{cfl,cflw}, we use all three formulae in \eqref{sda} to compute all the iterates for $G$, $H$ and $A$. For large-scale AREs, the corresponding SDA \cite{cw15,lclw13} compute only the iterates for $G$ and $H$ explicitly, leaving that for $A$ as a recursion. Interestingly, for the SDA\_h here for large-scale AREs with high-rank constant terms, only the iterates for $G$ are computed explicitly, leaving that for $H$ and $A$ as recursions. We may speculate that for large-scale AREs with high-rank terms $G$ and $H$, we compute {\it nothing} explicitly?}
%and we may control the growth of columns in $L^H_k(2)$. We choose not to, because the recursions for $A_k$ and $H_k$ can be applied efficiently enough with only compression on $B_k$, as seen in \eqref{a}--\eqref{mct}.
%Obviously, as the SDA converges, increasingly smaller but higher-rank components are added to $L^H_k(2)$. Apparent from (\ref{br}), the growth in the sizes and ranks of these iterates is potentially exponential. If the convergence is slow relative to the growth in $B_k$, the algorithm will fail.
Consider the QR decompositions with column pivoting of
\begin{equation}\label{ptcgh}
\begin{array}{rlc}
& L^{Gd}_{k}(2) = Q^G_{k} U^G_{k} + \widetilde{Q}^G_{k} \widetilde{U}^G_{k}, \ \ \
\| \widetilde{U}^G_{k} \| < \tau_g, \\
& L^{Hd}_{k}(2) = Q^H_{k} U^H_{k} + \widetilde{Q}^H_{k} \widetilde{U}^H_{k}, \ \ \
\| \widetilde{U}^H_{k} \| < \tau_h,
\end{array}
\end{equation}
where $\tau_g$, $\tau_h$ are some small tolerances controlling PTC of $L^{Gd}_{k}(2)$ and $L^{Hd}_{k}(2)$ respectively, $m^{g(2)}_{k}$ and $m^{h(2)}_{k}$ are the respective column numbers of $L^G_{k}(2)$ and $L^G_{k}(2)$ bounded above by some given $m_{\max}$.  Then their ranks satisfy
\[
r_{k}^{g} := \mathrm{rank} ( L^G_{k}(2) )\leq m^{g(2)}_{k} \leq m_{\max}, \ \
r_{k}^{h} := \mathrm{rank} ( L^H_{k}(2) )\leq m^{h(2)}_{k} \leq m_{\max}
\]
with $m_{\max}\ll N$. Also $Q^G_{k} \in \mathbb{R}^{N \times r_{k}^{g}}$ and $Q^H_{k} \in \mathbb{R}^{N \times r_{k}^{h}}$ are orthogonal
and $U^G_{k} \in \mathbb{R}^{r_{k}^{g} \times m^{hga}_{k-1}}$ and $U^H_{k} \in \mathbb{R}^{r_{k}^{h} \times m^{hga}_{k-1}}$ are full-rank and upper triangular with $m^{hga}_{k-1}= m^h_{k-1}+m^g_{k-1}+m^a_{k-1}$. Then $L^{Gd}_{k}$ and $L^{Hd}_{k}$ can be truncated and reorganized as
\begin{equation}\label{lghkp1}
\begin{array}{l}
L^{Gdt}_{k} =[ L^{Gd}_{k}(1), \  Q^G_{k}, \ L^{Gd}_{k}(3)]:=[ L^{Gdt}_{k}(1), \ L^{Gdt}_{k}(2), \ L^{Gdt}_{k}(3)]\in \mathbb{R}^{N \times m^g_k} , \\
L^{Hdt}_{k} =[ L^{Hd}_{k}(1), \ Q^H_{k},  \ L^{Hd}_{k}(3)]:=[ L^{Hdt}_{k}(1), \ L^{Hdt}_{k}(2), \ L^{Hdt}_{k}(3)]\in \mathbb{R}^{N \times m^h_k}
\end{array}
\end{equation}
with $m^g_k=r^g_k+km^a$ and $m^h_k=r^h_k+km^a$, respectively.

Similarly, recalling the delated forms in \eqref{lad1k} and \eqref{lad2k} in Appendix B, $L^{Ad}_{1, k}$ and $L^{Ad}_{2, k}$ will also be divided to two parts
\[
\begin{array}{l}
L^{Ad}_{1,k} = [L^{Ad}_{1,k}(1),  \ \ L^{Ad}_{1,k}(2)] \ \ \mbox{and} \ \
L^{Ad}_{2,k} = [L^{Ad}_{2,k}(1), \ \ L^{Ad}_{2,k}(2)]
\end{array}
\]
with
\[
\begin{array}{l}
L^{Ad}_{1,k}(1) = L^{Gd}_{k}(2),  \ \ L^{Ad}_{1,k}(2) = \Pi_{i=k-1}^0D_i^{A GH}L_{10}^A,\\
L^{Ad}_{2,k}(1) = L^{Hd}_{k}(2), \ \
L^{Ad}_{2,k}(2) = \Pi_{i=k-1}^0D_i^{A^\top HG}L_{20}^A.
\end{array}
\]
Since $L^{Gd}_{k}(2)$ and $L^{Hd}_{k}(2)$ have been compressed to $Q^{G}_{k}$ and $Q^{H}_{k}$, respectively. Then one has the truncated and compressed factors
\begin{equation}\label{lad12kp1}
\begin{array}{l}
L^{Adt}_{1,k} =[Q^{G}_{k}, L^{Ad}_{1,k}(2)]=[L^{Gdt}_{k}(2), L^{Adt}_{1,k}(2)]:=[L^{Adt}_{1,k}(1), L^{Adt}_{1,k}(2)]\in \mathbb{R}^{N \times m^{a_1}_k}, \\
L^{Adt}_{2,k} =[Q^{H}_{k}, L^{Ad}_{2,k}(2)]=[L^{Hdt}_{k}(2), L^{Adt}_{2,k}(2)]:=[L^{Adt}_{2,k}(1), L^{Adt}_{2,k}(2)]\in \mathbb{R}^{N \times m^{a_2}_k}
\end{array}
\end{equation}
with $m^{a_1}_k=r^g_k+m^a$ and $m^{a_2}_k=r^h_k+m^a$, finishing the PTC process for the low-rank factors in the $k$-th iteration.

It is worth noting that the above PTC process can proceed to the next iteration. In fact, one has
\[
\begin{array}{l}
L^{G}_{k+1} =[L_{k}^{Gdt}, \  L_{1,k}^{Adt},  \ D_{k}^{AGH}L_{k}^{Gdt}, \ D_{k}^{AGHG}L_{k}^{Hdt},  \ \ D_{k+1}^{AGHG}L_{2,k}^{Adt}], \\
L^{H}_{k+1} =[L_{k}^{Hdt}, \  L_{2,k}^{Adt},  \ D_{k}^{A^\top HG}L_{k}^{Hdt}, \ D_{k}^{A^\top HGH}L_{k}^{Gdt},  \ \ D_{k}^{A^\top HGH}L_{1,k}^{Adt}]
\end{array}
\]
after the $k$-th PTC.
As $L^{Adt}_{1,k}(1)$ is  equivalent to  $L^{Gdt}_{k}(2)$ and $L^{Adt}_{2,k}(1)$ is equivalent to  $L^{Hdt}_{k}(2)$, one can deflate $L^{G}_{k+1}$ and $L^{H}_{k+1}$ to
\[%\begin{equation}\label{lghdkp2}
\begin{array}{c}
L^{Gd}_{k+1}=[L^{Gd}_{k+1}(1), L^{Gd}_{k+1}(2), \ L^{Gd}_{k+1}(3)],  \ \ \
L^{Hd}_{k+1}=[L^{Hd}_{k+1}(1), L^{Hd}_{k+1}(2) L^{Hd}_{k+1}(3)]
\end{array}
\]%\end{equation}
with
\[
\begin{array}{c}
L^{Gd}_{k+1}(1)=[L_{k}^{Gdt}(1), \ L_{k}^{Gdt}(3)],\ \ L^{Hd}_{k+1}(1)=[L_{k}^{Hdt}(1),\ L_{k}^{Hdt}(3)], \\
L^{Gd}_{k+1}(2)=[L_{1,k}^{Adt},  \ D_{k}^{AGH}L_{k}^{Gdt}, \ D_{k}^{AGHG}L_{k}^{Hdt}],\ \
L^{Hd}_{k+1}(2)=[L_{2,k}^{Adt},  \ D_{k}^{A^\top HG}L_{k}^{Hdt}, \ D_{k}^{A^\top HGH}L_{k}^{Gdt}],\\
L^{Gd}_{k+1}(3)= D_{k}^{AGHG}L_{2,k}^{Adt}(2), \ \
L^{Hd}_{k+1}(3)= D_{k}^{A^\top HGH}L_{1,k}^{Adt}(2).
\end{array}
\]
Applying PTC to $L^{Gd}_{k+1}(2)$ and $L^{Hd}_{k+1}(2)$ respectively again, one has
\begin{equation}\label{lghdkp2}
\begin{array}{l}
L^{Gdt}_{k+1} =[L_{k+1}^{Gd}(1), \ Q_{k+1}^{G} \  L_{k+1}^{Gd}(3)]:=[L^{Gdt}_{k+1}(1), L^{Gdt}_{k+1}(2), \ L^{Gdt}_{k+1}(3)], \\
L^{Hdt}_{k+1} =[L_{k+1}^{Hd}(1),\ Q_{k+1}^{G}, \ L_{k+1}^{Hd}(3) ]:=[L^{Hdt}_{k+1}(1), L^{Hdt}_{k+1}(2), \ L^{Hdt}_{k+1}(3)],
\end{array}
\end{equation}
where
$Q^G_{k+1} \in \mathbb{R}^{N \times r_{k+1}^{g}}$ and $Q^H_{k+1} \in \mathbb{R}^{N \times r_{k+1}^{h}}$ are unitary matrices from QR decomposition and the PTC in the $(k+1)$-th iteration is completed.

\vspace{.2cm}
{\bf PTC for kernels.} Define orthogonal matrices
\[
\begin{array}{l}
\widehat{U}^A_{1,k} =  U^G_{k}\oplus I_{m^a}, \  \ \widehat{U}^G_{k} = I_{(k-1)m^a}\oplus U^G_{k}\oplus I_{m^a},\\
\widehat{U}^A_{2,k} =  U^H_{k}\oplus I_{m^a}, \ \
\widehat{U}^H_{k} = I_{(k-1)m^a}\oplus U^H_{k}\oplus I_{m^a},
\end{array}
\]
with $U^G_{k}$ and $ U^H_{k}$ in \eqref{ptcgh}. Then the truncated and compressed kernels are
\begin{equation}\label{kdtgha}
\begin{array}{l}
K^{Gdt}_{k} := \widehat{U}^G_{k} K^{Gd}_{k} (\widehat{U}^G_{k})^\top\in \mathbb{R}^{m_{k}^{g} \times m_{k}^{g}}, \\
K^{Hdt}_{k} := \widehat{U}^H_{k} K^{Hd}_{k} (\widehat{U}^H_{k})^\top\in \mathbb{R}^{m_{k}^{h} \times m_{k}^{h}}, \\
K^{Adt}_{k} := \widehat{U}^A_{1,k} K^{Hd}_{k} (\widehat{U}^A_{2,k})^\top \in \mathbb{R}^{m_{k}^{g} \times m_{k}^{h}},
\end{array}
\end{equation}
respectively.

To eliminate items less than $O(\tau_g)$ and $O(\tau_h)$ in low-rank factors and kernels, an additional monitoring step is imposed after PTC process.
Specifically, the last item $D^{AGHG}_{k-2}\Pi_{i=0}^{k-3}D_i^{A^\top GH}L_{20}^A$ in $L^{Gdt}_{k}$ (or $D^{A^\top HGH}_{k-2}\Pi_{i=0}^{k-3}D_i^{AGH}L_{10}^A$ in $L^{Gdt}_{k}$) will be discarded if its norm is less than $O(\tau_g)$ (or $O(\tau_h)$). Similarly,  $\Pi_{i=k-1}^{0}D_i^{AGH}L_{10}^A$ in $L^{Ad}_{1,k}(2)$ (or $\Pi_{i=k-1}^{0}D_i^{A^\top HG}L_{20}^A$ in $L^{Ad}_{2,k}(2)$) will be abandoned if its norm is less than $O(\tau_g)$ (or  $O(\tau_h)$).
In this way,  the growth of column dimension in low-rank factors $L^{Gdt}_{k}$, $L^{Hdt}_{k}$, $L^{Adt}_{1,k}$ and $L^{Adt}_{2,k}$, as well as kernels $K_{k}^{Gdt}$, $K_{k}^{Hdt}$, $K_{k}^{Adt}$, will be controlled efficiently with sacrificing a hopefully negligible bit of accuracy. Additionally, their sizes after PTC process will be further restricted by setting a reasonable upper bound $m_{\max}$.

\section{Algorithm and implementations}

\subsection{Computation of residual}
Define
\[
  \begin{array}{c}
    \widetilde{D}^{HG}_k = (I+D^H_kD^G_0)^{-1},\ \
    \widetilde{D}^{HGH}_k = \widetilde{D}^{HG}_kD^H_k,\ \
    \widetilde{D}^{GHG}_k = D^G_0\widetilde{D}^{HG}_k
  \end{array}
\]%end{equation}
and
\[
  \begin{array}{c}
    \widetilde{K}^{H}_k= (I+K^{H}_k(L^{H}_k)^\top\widetilde{D}^{GHG}_kL^{H}_k)^{-1}K^{H}_k.
  \end{array}
\]%end{equation}

With the current approximated solution $H_k = D^{H}_k + L^{H}_k K^{H}_k (L^{H}_k)^\top$, the residual for DARE \eqref{d} is
\[%begin{equation}
\begin{array}{rlc}
\mathcal{D} (H_k) &=
-H_k + A^\top  \Big( \widetilde{D}^{HGH}_k  + \widetilde{D}^{HG}_k L^H_k\widetilde{K}^{H}_k (\widetilde{D}^{HG}_k L^H_k)^\top \Big) A + H \\
&:= D^R_k + L^R_k K^R_k (L^R_k)^\top,
\end{array}
\]%end{equation}
where the banded part, the low-rank part and the kernel are
\[
D^R_k = D^H_0- D^H_k +(D^A_0)^\top D^H_k(I+D^G_0D^H_k)^{-1}D^A_0,\]
\[
L^R_k = [L^A_{20}, \ (D^A_0)^\top\widetilde{D}^{HGH}_k L^A_{10}, \ (D^A_0)^\top\widetilde{D}^{HG}_k L^H_k, \ L^H_k],
\]
\begin{equation}\label{krk}
K^R_k = \hspace{-.1cm}
\begin{array}{c@{\hspace{-2pt}}l}
\begin{array}{rr}
 \hspace{1cm} m^a   \hspace{1cm} m^a \hspace{.8cm}  m^{h}_k    \hspace{.8cm} m^{h}_k   \hspace{.2cm}&
\end{array}
&  \\
\left[
\begin{array}{cccc}
  \widetilde{K}^{A^\top HGHA}_k & I_{m^a}  &\widetilde{K}^{A^\top HG}_k &0 \\
 I_{m^a} & 0 & 0 & 0\\
 (\widetilde{K}^{A^\top HG}_k)^\top & 0 & \widetilde{K}^H_k & 0 \\
 0& 0 & 0 & -K^H_k
\end{array}
\right]
&
\begin{array}{l}
m^a \\ m^a \\ m^{h}_k \\ m^{h}_k \\
\end{array}
\end{array}
\end{equation}
respectively, and
\[
\widetilde{K}^{A^\top HG}_k= (L^A_{10})^\top \widetilde{D}^{HG}_k L^H_{k}\cdot \widetilde{K}^{H}_k,
\]
\[
\widetilde{K}^{A^\top HGHA}_k= (L^A_{10})^\top \widetilde{D}^{HGH}_k L^A_{10} + \widetilde{K}^{A^\top HG}_k\cdot \Big((L^A_{10})^\top \widetilde{D}^{HG}_k L^H_{k}\Big)^\top.
\]
It is not difficult to see that the main flops counts in the kernel $K^R_k$ lies in forming matrices
\begin{equation}\label{ron}
 (L^A_{10})^\top \widetilde{D}^{HGH}_k L^A_{10}, \ \ (L^A_{10})^\top \widetilde{D}^{HG}_k L^H_k, \ \  (L^H_k)^\top \widetilde{D}^{GHG}_k L^H_k.
\end{equation}
To avoid calculating them in each iteration, we can firstly set
\begin{equation}\label{bterm}
\mbox{B\_RRes} = \frac{ \|D^R_k\|}{\|D^H_0\|+ \|D^H_k\| +\|D^A_0\|^2 \|D^H_k\| /\|I+D^H_0D^H_k\|} \leq \epsilon_b
\end{equation}
as a pre-terminated condition of FSDA with $\epsilon_b$ the band tolerance. This is feasible as the residual of ${\mathcal D}(H_k)$ comes from two relatively independent parts, i.e. the banded part and the low-rank part. When the pre-termination \eqref{bterm} is satisfied,  matrices in \eqref{ron} are then constructed, followed by the deflation, truncation and compression of the low-rank factor $L^R_k$. Specifically, the columns $L^A_{20}(:,1:m^a)$ can be shift to columns $L^H_k(:,1:m^a)$ such that $L^R_k$ is deflated to $L^{Rd}_k$, i.e.
 \begin{eqnarray}
 L_{k}^R
&=&{
\begin{tikzpicture}[every node/.style={inner sep=0pt}]
\matrix (A) [
matrix of math nodes,
left delimiter ={[},
right delimiter = {]}
]
{L_{20}^A, & \ \ (D^A_0)^\top \widetilde{D}_k^{HGH}L_{10}^A, & \ \ (D^A_0)^\top \widetilde{D}^{HG}L_{k}^H, & \ \  L_k^{H} \  \\
};
\draw[->,black] (A-1-1.north) --+(+90:2mm) -| (A-1-4.north) node [pos=.35, above=.4em, text width=5.5cm,  font=\tiny] {$1 : m^a  \rightarrow m^{h}_k+2m^a+1 : m^{h}_k+3m^a$};
\end{tikzpicture}
} \nonumber\\
&\stackrel{d}{\rightarrow}&{
\begin{tikzpicture}[every node/.style={inner sep=0pt}]
\matrix (A) [
matrix of math nodes,
left delimiter ={[},
right delimiter = {]}
]
{ \ \ \ \ \ \  & \ \ (D^A_0)^\top \widetilde{D}_k^{HGH}L_{10}^A, & \ \ (D^A_0)^\top \widetilde{D}^{HG}L_{k}^H, & \ \  L_k^{H} \  \\
};
\end{tikzpicture}
} \nonumber\\
&:=& L_{k}^{Rd}.\nonumber
\end{eqnarray}
Let $\widehat{I}_{m^a} = [I_{m^a}, 0,...0]\in {\mathbb R}^{m^a\times m^{h}_k}$, $\widehat{K}^{A^\top HG}_k=[(\widetilde{K}^{A^\top HG}_k)^\top, 0,..., 0]\in {\mathbb R}^{m^{h}_k\times m^{h}_k}$. The kernel $K^R_k$ in \eqref{krk} is correspondingly deflated as
\[
K^R_k \stackrel{d}{\rightarrow}
\begin{array}{c@{\hspace{-2pt}}l}
\begin{array}{rr}
m^a \hspace{1.4cm}  m^{h}_k    \hspace{1.3cm} m^{h}_k   \hspace{.2cm}&
\end{array}
&  \\
\left[
\begin{array}{cccc}
   0&  0 & \widehat{I}_{m^a}  \\
  0 & \widetilde{K}^H_k & \widehat{K}^{A^\top HG}_k \\
 (\widehat{I}_{m^a})^\top& (\widehat{K}^{A^\top HG}_k)^\top &  \widehat{K}^{A^\top HGHA}_k
\end{array}
\right]
&
\begin{array}{l}
 m^a \\ m^{h}_k \\ m^{h}_k \\
\end{array}
\end{array} :=K^{Rd}_k,
\]
where all elements in $\widehat{K}^{A^\top HGHA}_k$ are same to those in $K^H_k$ except $\widehat{K}^{A^\top HGHA}_k(1:m^a, 1:m^a) = \widetilde{K}^{A^\top HGHA}_k- K^H_k(1:m^a, 1:m^a)$.

After deflation, the truncation and compression are applied to $L^{Rd}_k$ with QR decomposition
\[
\begin{array}{rlc}
 L^{Rd}_k = Q^R_{k} U^R_{k} + \widetilde{Q}^R_{k} \widetilde{U}^R_{k}, \ \ \
\| \widetilde{U}^R_{k} \| < \tau_r,
\end{array}
\]
where $\tau_r$ is the given tolerance, $Q^R_{k}\in {\mathbb R}^{n\times r^{r}_k}$ is unitary and $U^R_{k}\in {\mathbb R}^{r^{r}_k\times n_k}$ is full-rank and upper triangular. Then the terminated condition of the whole algorithm is
\begin{equation}\label{lrterm}
\mbox{LR\_RRes} = \frac{ \|U^R_k K^{Rd}_k (U^R_k)^\top\|} {\|U^R_k\|^2\|K^{Rd}_k\|} \leq \epsilon_l
\end{equation}
with $\epsilon_l$ the low-rank tolerance.

\subsection{Algorithm and operation counts}
The  process of deflation and PTC together with the computation of residual \eqref{bterm} and \eqref{lrterm} are summarized in the following FSDA algorithm.

\begin{table}[H]
\begin{center}
\begin{tabular}{l}\hline\hline%\label{alg1}
{\bf Algoritm FSDA.} Solve Riccati Equations with High-rank $G$ and $H$.  \\ \hline
\ \ Inputs: Banded matrices $D^A_0$, $D^G_0$, $D^H_0$, low-rank factors $L^A_{10}$, $L^A_{20}$, $L^G_0$, $L^H_0$ and the iterative\\
\hspace{1.2cm} \ \ tolerance $tol$, truncation tolerances $\tau_g$, $\tau_h$, $\tau_r$ and upper bound $m_{\max}$, band tolerance $\epsilon_b$  \\
\hspace{1.2cm} \ \ and low-rank tolerance $\epsilon_l$. \\
Outputs: Sparse banded matrix $D^H$, low-rank matrix $L^H$ and the kernel $K^H$ with the stabilizing  \\
\hspace{1.2cm} \ \ solution $X^\ast \approx D^H +L^HK^H(L^H)^\top$.  \\
\hline
\hspace{.1cm} 1.  \ \ Set $D^G_1$, $D^H_1$, $D^A_1$ in \eqref{sda1} with $D^G_0$, $D^H_0$, $D^A_0$; Compute low-rank factors $L^{G}_1$, $L^{H}_1$, $L^{A}_{11}$, $L^{A}_{21}$ \\
\hspace{.1cm} \ \  \ \ \ \ in \eqref{lagh1} and kernels $K^{G}_1$, $K^{H}_1$, $K^{A}_1$ in \eqref{kg1}-\eqref{ka1}.\\
\hspace{.1cm} 2.  \ \ For $k= 2, ..., $ until convergence, do \\
%%%%%%%%%%%%%
\hspace{.1cm} 3.\qquad \qquad Compute the banded matrices $D^G_k$, $D^H_k$, $D^A_k$ with iteration format \eqref{sdak}.  \\
\hspace{.1cm} 4.\qquad \qquad Form components \eqref{kghk}-\eqref{kghgk} and construct kernels $K^G_k$, $K^G_k$ and $K^G_k$ in \eqref{kgkp1}-\eqref{kakp1}.\\
\hspace{.1cm} 5.\qquad \qquad Deflate  kernels $K^G_k\stackrel{d}{\rightarrow}K^{Gd}_k$, $K^H_k\stackrel{d}{\rightarrow}K^{Hd}_k$ and $K^A_k\stackrel{d}{\rightarrow}K^{Ad}_k$ in a way of Fig. 3.3-3.4.\\
\hspace{.1cm} 6.\qquad \qquad Deflate the low-rank factors  $L^G_k\stackrel{d}{\rightarrow}L^{Gd}_k$, $L^H_k\stackrel{d}{\rightarrow}L^{Hd}_k$,  $L^A_{1,k}\stackrel{d}{\rightarrow}L^{Ad}_{1,k}$ and
$L^A_{2,k}\stackrel{d}{\rightarrow}L^{Ad}_{2,k}$ \\ \hspace{.1cm} \qquad \qquad \ \ \ in \eqref{lgdk}-\eqref{lad2k}.\\
%%%%%%%%%%%%%%%%%%%%
\hspace{.1cm} 7.\qquad \qquad Partially truncate and compress $L^{Gd}_k$ and $L^{Hd}_k$ in \eqref{ptcgh} with accuracy $\tau_g$, $\tau_h$. \\
\hspace{.1cm} 8.\qquad \qquad Construct compressed low-rank factors $L^{Gdt}_k$, $L^{Hdt}_k$, $L^{Adt}_{1,k}$ and $L^{Adt}_{2,k}$ in  \eqref{lghkp1}-\eqref{lad12kp1}.  \\
\hspace{.1cm} 9.\qquad \qquad Construct compressed kernels $K^{Gdt}_k$, $K^{Hdt}_k$ and $K^{Adt}_k$ in \eqref{kdtgha}.  \\
%%%%%%%%%%%%%%
\hspace{-.08cm} 10.\qquad \qquad Evaluate the residual of the banded part B\_RRes in \eqref{bterm}.  \\
\hspace{-.08cm} 11.\qquad \qquad If  B\_RRes $<tol$, compute the residual of low-rank part LR\_RRes in \eqref{lrterm}. \\
\hspace{-.08cm} 12.\qquad \qquad \qquad If  LR\_RRes $<tol$, break, end. \\
\hspace{-.08cm} 13.\qquad \qquad end;  \\
\hspace{-.08cm} 14.\qquad \qquad  $K^G_{k}:=K^{Gdt}_{k}$, $K^H_{k}:=K^{Hdt}_{k}$, $K^A_{k}:=K^{Adt}_{k}$. \\
\hspace{-.08cm} 15.\qquad \qquad  $L^G_k:=L^{Gdt}_k$, $L^H_k:=L^{Hdt}_k$, $L^A_{1,k}
:=L^{Adt}_{1,k}$, $L^A_{2,k}:=L^{Adt}_{2,k}$. \\
\hspace{-.08cm} 16.\qquad \qquad  $k:=k+1$; \\
\hspace{-.08cm} 17. \quad  End (for) \\
\hspace{-.08cm} 18. \quad  Output $D^H_k$, $L^H_k$ and $K^H_k$. \\
\hline\hline
 \end{tabular}
\end{center}
\end{table}

\begin{remark}
  1. Elements with the absolute value less than $tol$ in the banded matrices $D^G_k$, $D^G_k$ and $D^G_k$ will be eliminated at each iteration.

2. The deflation process just merges some rows and columns in kernels $K^G_k$, $K^H_k$ and $K^A_k$ according to the overlapped columns in low-rank factors $L^G_k$, $L^H_k$, $L^A_{1,k}$ and $L^A_{2,k}$, without costing any flops.

3. The PTC is only imposed on $L^{Gd}_k(2)$ and $L^{Hd}_k(2)$, with column numbers of  $L^{Gd}_k(1)$ and $L^{Hd}_k(1)$ increasing linearly about $k$ and those of $L^{Gd}_k(3)$ and $L^{Hd}_k(3)$ remaining invariant. As same before, elements in $L^{Gd}_k(1)$, $L^{Hd}_k(1)$, $L^{Gd}_k(3)$ and $L^{Hd}_k(3)$ with the absolute value less than $tol$ will be removed to slim the columns of low-rank factors as much as possible.
\end{remark}

To further analyze the complexity and the memory of FSDA, the bandwidth of $D_k^A$, $D_k^G$ and $D_k^H$ at each iteration
are assumed to be $b_k^a$, $b_k^g$ and $b_k^h$ ($b_k^a, b_k^g, b_k^h\ll N$), respectively.  We also set
$b_k^{hg}=\max\{b_k^{h}, b_k^{g}\}$,
$b_k^{hga}=\max\{b_k^{h}, b_k^{g}, b_k^{a}\}$,
$m_k^a=\max\{m_k^{a_1}, m_k^{a_2}\}$ and $m^{hga}_{k-1}:= m^h_{k-1}+m^g_{k-1}+m^a_{k-1}$ for the convenience of counting flops. The table in Appendix C lists the required flops and memories for different components in the $k$-th iteration of FSDA, where the estimations are roughly the upper bounds due to the truncation errors $\tau_g$, $\tau_h$ and $\tau_r$.

\begin{landscape}
\subsection*{Appendix A.}
\begin{center}

\end{center}
\begin{eqnarray}
L_2^G&=& [\ \ \ \ \ \ \ \  \ L_1^G \ \ \ \ \mid \ \  \ \ \ \ L_{11}^A \ \ \ \ \ \mid \ \ \ \ \ \ \ \ \ D_1^{AGH}L_1^G  \ \ \ \ \ \ \ \ \mid \ \ \ \ \ \ \ \ \ \ \ \ D_1^{AGHG}L_{1}^H \ \ \ \ \ \ \ \ \mid \ \ \ \ \ \ \ \  \ \ D_1^{AGHG}L_{21}^A\ \ \ \ \ \ \ \ \ ] \nonumber\\
&=&{\tiny
\begin{tikzpicture}[every node/.style={inner sep=0pt}]
\matrix (A) [
matrix of math nodes,
%left delimiter =(,
%right delimiter = )
left delimiter ={[},
right delimiter = {]}
]
{
L_{10}^A, & D_0^{AGHG}L_{20}^A, \mid &L_{10}^A, & D_0^{AGH}L_{10}^A, \mid & D_1^{AGH}L_{10}^A, &D_1^{AGH}D_0^{AGHG}L_{20}^A, \mid & D_1^{AGHG}L_{20}^A, & D_1^{AGHG}D_0^{A^\top GHG}L_{10}^A, \mid & D_1^{AGHG}L_{20}^A, & D_1^{AGHG}D_0^{A^\top HG}L_{20}^A\\
};
\draw[->,black] (A-1-1.north) --+(+90:2mm) -| (A-1-3.north);
\draw[->,black] (A-1-9.north) --+(+90:2mm) -| (A-1-7.north);
\end{tikzpicture}
} \nonumber\\
&\stackrel{d}{\rightarrow}&{\tiny
\begin{tikzpicture}[every node/.style={inner sep=0pt}]
\matrix (A) [
matrix of math nodes,
%left delimiter =(,
%right delimiter = )
left delimiter ={[},
right delimiter = {]}
]
{
\qquad \ \ & D_0^{AGHG}L_{20}^A,\mid & {\bf L_{10}^A}, & {\bf D_0^{AGH}L_{10}^A}, \mid &{\bf D_1^{AGH}L_{10}^A}, &{\bf D_1^{AGH}D_0^{AGHG}L_{20}^A}, \mid&{\bf D_1^{AGHG}L_{20}^A}, &{\bf D_1^{AGHG}D_0^{A^\top GHG}L_{10}^A}, \mid& \qquad \qquad \ \ & D_1^{AGHG}D_0^{A^\top HG}L_{20}^A\\
};
\end{tikzpicture}
} \nonumber\\
&:=& L_2^{Gd}.\nonumber\\
&\quad&\nonumber\\
L_{12}^A&=&[ \ \ \ \ \ \ L_{11}^A \ \ \ \ \ \mid \ \ \ \ \ \ \ \ \ D_1^{AGH}L_1^G \ \ \ \ \ \ \ \ \ \mid \ \ \ \ \ \ \ \ \ \ \ \ D_1^{AGHG}L_{1}^H \ \ \ \ \ \ \ \ \ \mid \ \ \ \ \ \ \ D_1^{AGHG}L_{21}^A\ \ \ \ \ \  ] \nonumber\\
&=&{\tiny
\begin{tikzpicture}[every node/.style={inner sep=0pt}]
\matrix (A) [
matrix of math nodes,
left delimiter ={[},
right delimiter = {]}
]
{
L_{10}^A, & D_0^{AGH}L_{10}^A, \mid & D_1^{AGH}L_{10}^A, &D_1^{AGH}D_0^{AGHG}L_{20}^A, \mid & D_1^{AGHG}L_{20}^A, & D_1^{AGHG}D_0^{A^\top GHG}L_{10}^A, \mid & D_1^{AGH}L_{10}^A, & D_1^{AGH}D_0^{AGH}L_{10}^A\\
};
\draw[->,black] (A-1-7.north) --+(+90:2mm) -| (A-1-3.north);
\end{tikzpicture}
} \nonumber\\
&\stackrel{d}{\rightarrow}&{\tiny
\begin{tikzpicture}[every node/.style={inner sep=0pt}]
\matrix (A) [
matrix of math nodes,
left delimiter ={[},
right delimiter = {]}
]
{
{\bf L_{10}^A}, & {\bf D_0^{AGH}L_{10}^A}, \mid &{\bf D_1^{AGH}L_{10}^A}, &{\bf D_1^{AGH}D_0^{AGHG}L_{20}^A}, \mid &{\bf D_1^{AGHG}L_{20}^A}, &{\bf D_1^{AGHG}D_0^{A^\top GHG}L_{10}^A} \mid & \qquad \qquad \ \ & D_1^{AGH}D_0^{AGH}L_{10}^A\\
};
\end{tikzpicture}
} \nonumber\\
&:=& L_{12}^{Ad}.\nonumber\\
&\quad&\nonumber\\
L_2^H&=& [\ \ \ \ \ \ \ \  \ L_1^H \ \ \ \ \ \ \mid \ \  \ \ \ \ L_{21}^A \ \ \ \ \ \ \mid \ \ \ \ \ \ \ \ \ \ \ D_1^{A^\top HG}L_1^H \ \ \ \ \ \ \ \ \ \ \mid \ \ \ \ \ \ \ \ \ \ \ \ D_1^{A^\top HGH}L_{1}^G \ \ \ \ \ \ \  \ \ \mid \ \ \ \ \ \ \ \  \ \ D_1^{A^\top HGH}L_{11}^A\ \ \ \ \ \ \ \ \ ] \nonumber\\
&=&{\tiny
\begin{tikzpicture}[every node/.style={inner sep=0pt}]
\matrix (A) [
matrix of math nodes,
left delimiter ={[},
right delimiter = {]}
]
{
L_{20}^A, & D_0^{A^\top HGH}L_{10}^A, \mid &L_{20}^A, & D_0^{A^\top HG}L_{20}^A, \mid & D_1^{A^\top HG}L_{20}^A, &D_1^{A^\top HG}D_0^{A^\top HGH}L_{10}^A, \mid & D_1^{A^\top HGH}L_{10}^A, & D_1^{A^\top HGH}D_0^{AGHG}L_{20}^A, \mid & D_1^{A^\top HGH}L_{10}^A, & D_1^{A^\top HGH}D_0^{AGH}L_{10}^A\\
};
\draw[->,black] (A-1-1.north) --+(+90:2mm) -| (A-1-3.north);
\draw[->,black] (A-1-9.north) --+(+90:2mm) -| (A-1-7.north);
\end{tikzpicture}
} \nonumber\\
&\stackrel{d}{\rightarrow}&{\tiny
\begin{tikzpicture}[every node/.style={inner sep=0pt}]
\matrix (A) [
matrix of math nodes,
%left delimiter =(,
%right delimiter = )
left delimiter ={[},
right delimiter = {]}
]
{
\qquad \ \ & D_0^{A^\top HGH}L_{10}^A, \mid& {\bf L_{20}^A}, & {\bf D_0^{A^\top HG}L_{20}^A}, \mid &{\bf D_1^{A^\top HG}L_{20}^A}, &{\bf D_1^{A^\top HG}D_0^{A^\top HGH}L_{10}^A},\mid &{\bf D_1^{A^\top HGH}L_{10}^A}, & {\bf D_1^{A^\top HGH}D_0^{AGHG}L_{20}^A}, \mid & \qquad \qquad \qquad &  D_1^{A^\top HGH}D_0^{AGH}L_{10}^A\\
};
\end{tikzpicture}
} \nonumber\\
&:=& L_2^{Hd}.\nonumber\\
&\quad& \nonumber\\
L_{22}^A&=&[\ \ \ \ \ \ \ \ L_{21}^A \ \ \ \ \ \mid \ \ \ \ \ \ \ \ \ \ \ D_1^{A^\top HG}L_1^H \ \ \ \ \ \ \ \ \ \ \mid \ \ \ \ \ \ \ \ \ \ \ \ D_1^{A^\top HGH}L_{1}^G \ \ \ \ \ \ \ \ \mid \ \ \ \ \ \ \ \  \ \ D_1^{A^\top HG}L_{21}^A\ \ \ \ \ \ \ \ \ ] \nonumber\\
&=&{\tiny
\begin{tikzpicture}[every node/.style={inner sep=0pt}]
\matrix (A) [
matrix of math nodes,
left delimiter ={[},
right delimiter = {]}
]
{
{L_{20}^A}, & { D_0^{A^\top HG}L_{20}^A}, \mid&{ D_1^{A^\top HG}L_{20}^A}, &{ D_1^{A^\top HG}D_0^{A^\top HGH}L_{10}^A},\mid &{D_1^{A^\top HGH}L_{10}^A}, & { D_1^{A^\top HGH}D_0^{AGHG}L_{20}^A},\mid & D_1^{A^\top HG}L_{20}^A, & D_1^{A^\top HG}D_0^{A^\top HG}L_{10}^A\\
};
\draw[->,black] (A-1-7.north) --+(+90:2mm) -| (A-1-3.north);
\end{tikzpicture}
} \nonumber\\
&\stackrel{d}{\rightarrow}&{\tiny
\begin{tikzpicture}[every node/.style={inner sep=0pt}]
\matrix (A) [
matrix of math nodes,
%left delimiter =(,
%right delimiter = )
left delimiter ={[},
right delimiter = {]}
]
{
{\bf L_{20}^A}, & {\bf D_0^{A^\top HG}L_{20}^A}, \mid&{\bf D_1^{A^\top HG}L_{20}^A}, &{\bf D_1^{A^\top HG}D_0^{A^\top HGH}L_{10}^A}, \mid&{\bf D_1^{A^\top HGH}L_{10}^A}, & {\bf D_1^{A^\top HGH}D_0^{AGHG}L_{20}^A}, \mid& \qquad \qquad \ \ & D_1^{AGH}D_0^{AGH}L_{10}^A\\
};
\end{tikzpicture}
} \nonumber\\
&:=& L_{22}^{Ad}.\nonumber
\end{eqnarray}
\end{landscape}

%%%%%%%%%%%%%%%%%%%%%%%%%%%%%%%%%%%%%%%%%%%%
\begin{landscape}
\subsection*{Appendix B.}
\begin{eqnarray}
L_{k}^G
&=&{
\begin{tikzpicture}[every node/.style={inner sep=0pt}]
\matrix (A) [
matrix of math nodes,
left delimiter ={[},
right delimiter = {]}
]
{
\hspace{3.5cm}  L_{k-1}^G \hspace{3.5cm} \mid &  L_{1,k-1}^A \mid  & D_{k-1}^{AGH}L_{k-1}^G \mid & D_{k-1}^{AGHG}L_{k-1}^H  \mid & \ \ \ D_{k-1}^{AGHG}L_{2,k-1}^A  \ \ \ \\
};
\draw[->,black] (A-1-1.north) --+(+90:2mm) -| (A-1-2.north) node [pos=0.26, above=.5em, text width=5cm, font=\tiny] {$(k-1)m^a+1 : m^g_{k-1}-m^a \rightarrow 1 : m^g_{k-1}-km^a$};
\draw[->,black] (A-1-5.north) --+(+90:2mm) -| (A-1-4.north) node [pos=.25, above=.6em, text width=6cm,  font=\tiny] {$1 : m^{a_2}_{k-1}-m^a  \rightarrow m^h_{k-1}-m^{a_2}_{k-1}+1 : m^h_{k-1}-m^a$};
\end{tikzpicture}
} \nonumber\\
&\stackrel{d}{\rightarrow}&{\tiny
\begin{tikzpicture}[every node/.style={inner sep=0pt}]
\matrix (A) [
matrix of math nodes,
left delimiter ={[},
right delimiter = {]}
]
{
D_0^{AGHG}L_{20}^A, \ D_1^{AGHG}D_0^{A^\top GH}L_{20}^A, \ ..., \ D_{k-2}^{AGHG}\Pi_{i=0}^{k-3}D_i^{A^\top GH}L_{20}^A\ \mid & \ \  { L_{1,k-1}^A} \ \ \mid & \ \ { D_{k-1}^{AGH}L_{k-1}^G} \ \ \mid & \ \ \ { D_{k-1}^{AGHG}L_{k-1}^H} \ \ \mid & D_{k-1}^{AGHG}\Pi_{i=k-2}^{0}D_i^{A^\top GH}L_{20}^A\\
};
\end{tikzpicture}
} \nonumber\\
&:=& L_{k}^{Gd},\label{lgdk}\\
&\quad&\nonumber\\
L_{1,k}^A
&=&{
\begin{tikzpicture}[every node/.style={inner sep=0pt}]
\matrix (A) [
matrix of math nodes,
left delimiter ={[},
right delimiter = {]}
]
{
L_{1,k-1}^A \mid & \ \ D_{k-1}^{AGH}L_{k-1}^G \mid & \ \ D_{k-1}^{AGHG}L_{k-1}^H \mid   & \ \  D_{k-1}^{AGH}L_{1,k-1}^A \ \ \ \  \\
};
\draw[->,black] (A-1-4.north) --+(+90:2mm) -| (A-1-2.north) node [pos=.25, above=.6em, text width=5.5cm,  font=\tiny] {$1 : m^{a_1}_{k-1}-m^a  \rightarrow m^{g}_{k-1}-m^{a_1}_{k-1}+1 : m^{g}_{k-1}-m^{a}$};
\end{tikzpicture}
} \nonumber\\
&\stackrel{d}{\rightarrow}&{
\begin{tikzpicture}[every node/.style={inner sep=0pt}]
\matrix (A) [
matrix of math nodes,
left delimiter ={[},
right delimiter = {]}
]
{
{ L_{1,k-1}^A} \mid & \ \ { D_{k-1}^{AGH}L_{k-1}^G} \mid & \  { D_{k-1}^{AGHG}L_{k-1}^H} \ \mid & \ \ \Pi_{i=k-1}^0D_i^{AGH}L_{10}^A\\
};
\end{tikzpicture}
} \nonumber\\
&:=& L_{1,k}^{Ad},\label{lad1k}\\
&\quad&\nonumber\\
L_{k}^H
&=&{
\begin{tikzpicture}[every node/.style={inner sep=0pt}]
\matrix (A) [
matrix of math nodes,
left delimiter ={[},
right delimiter = {]}
]
{
\hspace{3.5cm}  L_{k-1}^H \hspace{3.5cm}  \mid&   \ L_{2,k-1}^A  \mid  & D_{k-1}^{A^\top HG}L_{k-1}^H \mid & D_{k-1}^{A^\top HGH}L_{k-1}^G \mid   & \ \ D_{k-1}^{A^\top HGH}L_{1,k-1}^A  \ \ \ \\
};
\draw[->,black] (A-1-1.north) --+(+90:2mm) -| (A-1-2.north) node [pos=0.26, above=.5em, text width=5cm, font=\tiny] {$(k-1)m^a+1 : m^h_{k-1}-m^a \rightarrow 1 : m^h_{k-1}-km^a$};
\draw[->,black] (A-1-5.north) --+(+90:2mm) -| (A-1-4.north) node [pos=.25, above=.6em, text width=6cm,  font=\tiny] {$1 : m^{a_1}_{k-1}-m^a  \rightarrow m^g_{k-1}-m^{a_1}_{k-1}+1 : m^g_{k-1}-m^a$};
\end{tikzpicture}
} \nonumber\\
&\stackrel{d}{\rightarrow}&{\tiny
\begin{tikzpicture}[every node/.style={inner sep=0pt}]
\matrix (A) [
matrix of math nodes,
left delimiter ={[},
right delimiter = {]}
]
{
D_0^{A^\top HGH}L_{10}^A, \ D_1^{A^\top HGH}D_0^{AGH}L_{10}^A,..., D_{k-2}^{A^\top HGH}\Pi_{i=0}^{k-3}D_i^{AGH}L_{10}^A \mid &\ \ \ { L_{2,k-1}^A}\ \ \mid & \ \ { D_{k-1}^{A^\top HG}L_{k-1}^H} \ \ \mid & \ \ \ { D_{k-1}^{A^\top HGH}L_{k-1}^G} \ \ \ \mid & \  D_{k-1}^{A^\top HGH}\Pi_{i=k-2}^{0}D_i^{AGH}L_{10}^A\\
};
\end{tikzpicture}
} \nonumber\\
&:=& L_{k}^{Hd},\label{lhdk}\\
&\quad& \nonumber\\
L_{2,k}^A
&=&{
\begin{tikzpicture}[every node/.style={inner sep=0pt}]
\matrix (A) [
matrix of math nodes,
left delimiter ={[},
right delimiter = {]}
]
{L_{2,k-1}^A \mid & \ \ D_{k-1}^{A^\top HG}L_{k-1}^H \mid & \ \ D_{k-1}^{A^\top HGH}L_{k-1}^G \mid   & \ \  D_{k-1}^{A^\top HG}L_{2,k-1}^A \  \\
};
\draw[->,black] (A-1-4.north) --+(+90:2mm) -| (A-1-2.north) node [pos=.25, above=.6em, text width=5.5cm,  font=\tiny] {$1 : m^{a_2}_{k-1}-m^a  \rightarrow m^h_{k-1}-m^{a_2}_{k-1}+1 : m^h_{k-1}-m^a$};
\end{tikzpicture}
} \nonumber\\
&\stackrel{d}{\rightarrow}&{
\begin{tikzpicture}[every node/.style={inner sep=0pt}]
\matrix (A) [
matrix of math nodes,
left delimiter ={[},
right delimiter = {]}
]
{
{ L_{2,k-1}^A} \mid & \ \ { D_{k-1}^{A^\top HG}L_{k-1}^H} \mid & \ \ { D_{k-1}^{A^\top HGH}L_{k-1}^G} \mid & \ \ \Pi_{i=k-1}^0D_i^{A^\top HG}L_{20}^A  \\
};
\end{tikzpicture}
} \nonumber\\
&:=& L_{2,k}^{Ad}.\label{lad2k}
\end{eqnarray}
\end{landscape}

\begin{landscape}
\subsection*{Appendix C. {\rm Complexity and memory at $k$-th iteration in FSDA}}
\begin{center}
$${\footnotesize
\begin{tabular}{c c c} \hline
 Items & Flops & Memory  \\ \hline
& {Banded part} & \\   \hline
 %%%%%%%%%%%%%%%%%%%%%%%%%%
 $D_k^{AGH}$, $D_k^{{A^\top HG}^\ast}$ & { $4N(2b^{hg}_{k-1}+1)^2+b^{hg}_{k-1}b^{a}_{k-1})$}  &{ $2N(2b^{hga}_{k-1}+1)$} \\
 %%%%%%%%%%%%%%%%%%%%%%%%%%
 $D_k^{G}$, $D_k^{H}$, $D_k^{A}$ & { $4N(2b^{g}_{k-1}+1)(2b^{hga}_{k-1}+1)$}  &{$2N(2b^{hga}_{k-1}+1)$} \\ \hline
 & Low-rank part and kernels & \\   \hline
 %%%%%%%%%%%%%%%%%%%%%%%%%%
 $D_{k-1}^{AGH}L_{k-1}^{G}$, $D_{k-1}^{AGHG}L_{k-1}^{H}$, $D_{k-1}^{AGHG}L_{2,k-1}^{A}$  & $2Nb_{k-1}^{hga}(m_{k-1}^g+m_{k-1}^h+m_{k-1}^a)$ & {$(m_{k-1}^g+m_{k-1}^h+m_{k-1}^a)N$} \\
 %%%%%%%%%%%%%%%%%%%%%%%%%%%%%%%%%%%
 $D_{k-1}^{A^\top HG}L_{k-1}^{H}$, $D_{k-1}^{A^\top HGH}L_{k-1}^{G}$, $D_{k-1}^{A^\top HGH}L_{1,k-1}^{A}$ & $2Nb_{k-1}^{hga}(m_{k-1}^g+m_{k-1}^h+m_{k-1}^a)$ & {$(m_{k-1}^g+m_{k-1}^h+m_{k-1}^a)N$}\\
 %%%%%%%%%%%%%%%%%%%%%%%%%%%%%%%%%%%%%%
 $\Theta_{k-1}^H$, $\Theta_{k-1}^G$, $\Theta_{k-1}^{HG}$ & $\begin{array}{c}
   2N(b_{k-1}^{hg}(m_{k-1}^h+m_{k-1}^g)+b_{k-1}^{hg}m_{k-1}^g\\
   + (m_{k-1}^h)^2+(m_{k-1}^g)^2+m_{k-1}^gm_{k-1}^h)
 \end{array}$  &{ $(m_{k-1}^h)^2+(m_{k-1}^g)^2+m_{k-1}^hm_{k-1}^g$} \\
 %%%%%%%%%%%%%%%%%%%%%%%%%%%%%%
 $\Theta_{k-1}^{A}$, $\Theta_{1,k-1}^{A}$, $\Theta_{2,k-1}^{A}$ &{ $2N(2b_{k-1}^{hg}m_{k-1}^a+b_{k-1}^{hg}m_{k-1}^a +3(m_{k-1}^a)^2)$ } & {$3(m_{k-1}^a)^2$} \\
 %%%%%%%%%%%%%%%%%%%%%%%%%%%%%%
 $\Theta_{1,k-1}^{AH}$, $\Theta_{1,k-1}^{AG}$   &{ $2N(b_{k-1}^{hg}(m_{k-1}^h+m_{k-1}^g) +m_{k-1}^a(m_{k-1}^h+m_{k-1}^g))$} & { $m_{k-1}^a(m_{k-1}^h+m_{k-1}^g)$ } \\
 %%%%%%%%%%%%%%%%%%%%%%%%%%%%%%
 $\Theta_{2,k-1}^{AH}$, $\Theta_{2,k-1}^{AG}$  &{ $2N(b_{k-1}^{hg}(m_{k-1}^h+m_{k-1}^g) +m_{k-1}^a(m_{k-1}^h+m_{k-1}^g))$} &{ $m_{k-1}^a(m_{k-1}^h+m_{k-1}^g)$}  \\
 %%%%%%%%%%%%%%%%%%%%%%%%%%%%%%%
 $K_{k-1}^{AGHG}$  &{ $(m_{k-1}^a)^2(m_{k-1}^h+m_{k-1}^g)+
 m_{k-1}^a(m_{k-1}^h+m_{k-1}^g)^2$}&
 { $m_{k-1}^a(m_{k-1}^h+m_{k-1}^g)$}  \\
 %%%%%%%%%%%%%%%%%%%%%%%%%%%%%%%
 $K_{k-1}^{AGHGA^\top}$\hspace{-.2cm}, $K_{k-1}^{A^\top HGHA}$, $K_{k-1}^{AGHA}$  &{ $6(m_{k-1}^a)^2(2m_{k-1}^a+m_{k-1}^h+m_{k-1}^g)$}&{ $3(m^a_{k-1})^2$}  \\
  %%%%%%%%%%%%%%%%%%%%%%%%%%%%%%
 $K_{k-1}^{A^\top HGH}$  &{ $2(m_{k-1}^a)(m_{k-1}^a+m_{k-1}^h)
 (m_{k-1}^a+m_{k-1}^h+m_{k-1}^g)$}&
 { $m_{k-1}^a(m_{k-1}^h+m_{k-1}^g)$}  \\
  %%%%%%%%%%%%%%%%%%%%%%%%%%%%%%
 $K_{k-1}^{AGH}$, $K_{k-1}^{A^\top GH}$ &{ $2(m_{k-1}^a)(m_{k-1}^a+m_{k-1}^h)^2$ }& {$2m_{k-1}^a(m_{k-1}^h+m_{k-1}^g)$} \\
  %%%%%%%%%%%%%%%%%%%%%%%%%%%%%%%%%%
 $K_{k-1}^{GH}$, $K_{k-1}^{GHG}$, $K_{k-1}^{HGH}$ &{ $8(m_{k-1}^h+m_{k-1}^g)^3/3$ }& {$3(m_{k-1}^h+m_{k-1}^g)^2$} \\
 %%%%%%%%%%%%%%%%%%%%%%%%%%%%%%%%%%
 $Q_{k}^G$, $Q_{k}^{H^{\ast\ast}}$ &{ $4(m_{k-1}^a+m_{k-1}^g+m_{k-1}^h)^2(N-m_{k-1}^a+m_{k-1}^g+m_{k-1}^h)$ }& {$(r_{k}^{h}+r_{k}^{g})N$}\\
 %%%%%%%%%%%%%%%%%%%%%%%%%%%%%%%%%%%%%
 $U_{k}^G$, $U_{k}^H$,  &{ $4(m_{k-1}^a+m_{k-1}^g+m_{k-1}^h)r^{g}_{k-1}
 (N-m_{k-1}^a+m_{k-1}^g+m_{k-1}^h)$ }& {$(r_{k}^{g}+r_{k}^{h})\times m^{hga}_{k-1}$}\\
 %%%%%%%%%%%%%%%%%%%%%%%%%%%%%%%%%%%%%
 $K_{k}^{Gdt}$ &{ $12(m_{k-1}^a+m_{k-1}^g+m_{k-1}^h)^2 r^g_{k-1}$} & {$(m_{k}^{g})^2$}\\
 %%%%%%%%%%%%%%%%%%%%%%%%%%%%%%%%%%%%%
 $K_{k}^{Hdt}$, &{ $12(m_{k-1}^a+m_{k-1}^g+m_{k-1}^h)^2r^h_{k-1}$}& {$(m_{k}^{h})^2$}\\
 %%%%%%%%%%%%%%%%%%%%%%%%%%%%%%%%%%%%%
 $K_{k}^{Adt}$ &{ $6(m_{k-1}^a+m_{k-1}^g+m_{k-1}^h)^2(r^g_{k-1}+r^g_{k-1})$}& {$m_{k}^{g}m_{k}^{h}$}\\ \hline
 %%%%%%%%%%%%%%%%%%%%%%%%%%%%%%%%%%%%%
 & Residual part & \\   \hline
 %%%%%%%%%%%%%%%%%%%%%%%%%%%%%%%%%%%%
 \vspace{.1cm}
 $(D^A_0)^\top \widetilde{D}_{k}^{HGH}L^A_{10}$, $(D^A_0)^\top \widetilde{D}_{k}^{HG}L^H_{k}$ & { $2b^{hg}_{k}(m^a+m^h_k)N$} &{$(m^h_k+m^a)N$ } \\
 %%%%%%%%%%%%%%%%%%%%%%%%%%%%%%%%%%
 $(L^H_k)^\top \widetilde{D}_{k}^{GHG}L^H_{k}$& { $2b^{hg}_{k}(m^a+m^h_k)N$} &{$(m^h_k)^2$ } \\
 %%%%%%%%%%%%%%%%%%%%%%%%%%%%%%%%%%
 \vspace{.05cm}
 $\widetilde{K}_{k}^{H^\ast}$& { $8(m^h_k)^2/3$} &{$(m^h_k)^2$ } \\
 %%%%%%%%%%%%%%%%%%%%%%%%%%%%%%%%%%
 \vspace{.05cm}
 $\widetilde{K}_{k}^{A^\top HG}$& { $2b^{hg}_{k}(m^a+m^h_k)N$} &{$m^am^h_k$ } \\
 %%%%%%%%%%%%%%%%%%%%%%%%%%%%%%%%%%
 $\widetilde{K}_{k}^{A^\top HGHA}$& {$2m^a(b^{hg}_{k}+m^a)N + 2(m^a)^2m^h_k$} &{$(m^a)^2$ } \\
 %%%%%%%%%%%%%%%%%%%%%%%%%%%%%%%%%%
 \vspace{.05cm}
 $Q_{k}^{R^{\ast\ast}}$ &{ $2(m^a+2m_{k}^h)^2(N-m^a-2m_{k}^h)$ }& {$r_k^rN$}\\
 %%%%%%%%%%%%%%%%%%%%%%%%%%%%%%%%%%%%%
 \vspace{.05cm}
 $U_{k}^R$ &{ $2(m^a+2m_{k}^h)r_k^r(N-m^a-2m_{k}^h)$ }& {$r_k^r(m^a+2m_{k}^h)$}\\
 %%%%%%%%%%%%%%%%%%%%%%%%%%%%%%%%%%%%%
 $U_{k}^{R}K_{k}^{Rd}(U_{k}^{R})^\top$ &{ $2(m^a+2m_{k}^h)r_k^r(r_k^r+m^a+2m_{k}^h)$} & {$(r_k^r)^2$}\\
 %%%%%%%%%%%%%%%%%%%%%%%%%%%%%%%%%%%%%
 \hline\hline
 %%%%%%%%%%%%%%%%%%%%%%%%
\end{tabular}
}$$
\end{center}
{\footnotesize ${}^\ast$ LU factorization and Gaussian elimination is used \cite{ag9411}.

\noindent ${}^{\ast\ast}$ Householder QR decomposition is used \cite{m}.}

\end{landscape}

\end{document}